\documentclass[12pt]{article}
\usepackage{microtype} \DisableLigatures{encoding = *, family = * }
\usepackage{subeqn}
\usepackage{graphicx}
\usepackage{ifpdf}
\ifpdf  \fi
\usepackage{amsfonts,amssymb,slashbox}
\usepackage{mathptmx,helvet,courier,makeidx,multicol,footmisc}
\usepackage[numbers]{natbib}
\bibpunct{[}{]}{;}{n}{,}{,}

\newcommand{\commentout}[1]{}
\newcommand{\ba}{\begin{array}}
        \newcommand{\ea}{\end{array}}
\newcommand{\bc}{\begin{center}}
        \newcommand{\ec}{\end{center}}
\newcommand{\bdm}{\begin{displaymath}}
        \newcommand{\edm}{\end{displaymath}}
\newcommand{\bds} {\begin{description}}
        \newcommand{\eds} {\end{description}}
\newcommand{\ben}{\begin{enumerate}}
        \newcommand{\een}{\end{enumerate}}
\newcommand{\beq}{\begin{equation}}
        \newcommand{\eeq}{\end{equation}}
\newcommand{\bfg} {\begin{figure}[h]}
        \newcommand{\efg} {\end{figure}}
\newcommand{\bi} {\begin {itemize}}
        \newcommand{\ei} {\end {itemize}}
\newcommand{\bqn}{\begin{eqnarray}}
        \newcommand{\eqn}{\end{eqnarray}}
\newcommand{\bqs}{\begin{eqnarray*}}
        \newcommand{\eqs}{\end{eqnarray*}}
\newcommand{\bsl} {\begin{slide}[8.8in,6.7in]}
        \newcommand{\esl} {\end{slide}}
\newcommand{\bsq}{\begin{subequations}}
        \newcommand{\esq}{\end{subequations}}
\newcommand{\bss} {\begin{slide*}[9.3in,6.7in]}
        \newcommand{\ess} {\end{slide*}}
\newcommand{\btb} {\begin {table}}
        \newcommand{\etb} {\end {table}}
\newcommand{\bvb}{\begin{verbatim}}
        \newcommand{\evb}{\end{verbatim}}

\newcommand{\m}{\mbox}

\newcommand {\pd}[2] {{\frac {\partial {#1}} {\partial {#2}}}}

\newcommand{\cas}[1]{{{\left \{ \ba #1 \ea \right. }}}

\newcommand{\reff}[1] {{{Figure \ref {#1}}}}
\newcommand{\refe}[1] {{{(\ref {#1})}}}
\newcommand{\reft}[1] {{{\textbf{Table} \ref {#1}}}}
\def\la      {{\lambda}}

\def\pmb#1{\setbox0=\hbox{$#1$}%
   \kern-.025em\copy0\kern-\wd0
   \kern.05em\copy0\kern-\wd0
   \kern-.025em\raise.0433em\box0 }

\def\eop{{\hfill $\blacksquare$}}

\def\Po  {{Poincar\'e }}

\newtheorem{theorem}{Theorem}[section]

\newtheorem{lemma}[theorem]{Lemma}
\newtheorem{corollary}[theorem]{Corollary}

\usepackage{ulem}
\begin {document}

\title{Stability and bifurcation in network traffic flow: A Poincar\'e map approach} 
\author{Wen-Long Jin \footnote{Department of Civil and Environmental Engineering, California Institute for Telecommunications and Information Technology, Institute of Transportation Studies, 4000 Anteater Instruction and Research Bldg, University of California, Irvine, CA 92697-3600. Tel: 949-824-1672. Fax: 949-824-8385. Email: wjin@uci.edu. Corresponding author}}
\maketitle

\begin{abstract}
Previous studies have shown that, in a diverge-merge network with two intermediate links (the DM network), the kinematic wave model always admits stationary solutions under constant boundary conditions, but periodic oscillations can develop from empty initial conditions. Such contradictory observations suggest that the stationary states be unstable.
In this study we develop a new approach to investigate the stability property of traffic flow in this and other networks.
Based on the observation that kinematic waves propagate in a circular path when only one of the two intermediate links is congested, we  derive a one-dimensional, discrete \Po map in the out-flux at a \Po section. We then prove that the fixed points of the \Po map correspond to stationary flow-rates on the two links. With Lyapunov's first method, we demonstrate that the \Po map can be finite-time stable, asymptotically stable, or unstable. When unstable, the map is found to have periodical points of period two, but no chaotic solutions. Comparing the results with those in existing studies, we conclude that the \Po map can be used to represent network-wide dynamics in the kinematic wave model. We further analyze the bifurcation in the stability of the \Po map caused by varying route choice proportions.  We further apply the \Po map approach to analyzing traffic patterns in more general $(DM)^n$ and beltway networks, which are sufficient and necessary structures for network-induced unstable traffic and gridlock, respectively. This study demonstrates that the \Po map approach can be efficiently applied to analyze traffic dynamics in any road networks with circular information propagation and provides new insights into unstable traffic dynamics caused by interactions among network bottlenecks.

\end{abstract}


\section{Introduction}
Observations have demonstrated that stop-and-go patterns frequently arise in congested traffic \citep{mauch2002freeway}. As large speed variations in such oscillatory traffic patterns are challenging for human drivers, accident likelihood \citep{oh2001real} and vehicle emissions \citep{barth2008real} can be significantly increased. Therefore, the properties, mechanism, and control of stop-and-go traffic have been subject to many studies in the transportation engineering field. For example, methods have been proposed to measure the magnitudes, periods, and propagation of traffic oscillations \citep{li2010measurement}. It was shown that unstable car-following behaviors can lead to phantom jams and stop-and-go traffic in a platoon of vehicles  \citep[e.g.][]{herman1959traffic,kerner1994cluster,sugiyama2008traffic}. Lane-changing, merging, and diverging activities are also significant contributors of the formation and growth of traffic oscillations \citep{ahn2007freeway,ahn2010merging}. Generally, such instability is caused by individual vehicles' local behaviors.  

With many metropolitan areas becoming oversaturated, such network bottlenecks as merges and diverges play a more important role in shaping traffic dynamics \citep{daganzo1999remarks,daganzo1999phase}, and complicated network-wide traffic dynamics can be caused by interactions among different network bottlenecks. In \citep{daganzo1996gridlock}, a beltway network was shown to become totally gridlocked with certain merging priorities and diverging ratios. In \citep{daganzo2011bifurcations}, it was shown that a double-ring network can have multiple stationary states at the same density, and bifurcations are shown to exist in the macroscopic fundamental diagram. 
A series of studies have demonstrated that periodical oscillations with a period of about 10 minutes can occur in a diverge-merge network with two intermediate links in \reff{fig:dm2}: in \citep[][Section 7.3]{jin2003dissertation},  damped and persistent oscillatory traffic patterns were first observed with a commodity-based Cell Transmission Model (CTM), even when the network is initially empty and has constant demand patterns and route choice proportions; in \citep{jin2005paramics}, such network-wide oscillations were replicated with a microscopic traffic simulator Paramics, and it suggests that such periodic oscillations are intrinsic properties of the network, not caused by the specific traffic flow models on a link or at the junctions; in \citep{jin2009network}, the mechanism and network conditions for the occurrence of persistent periodic oscillations were identified with detailed analyses of shock and rarefaction waves in the network; in \citep{jin2012statics}, it was shown that the DM network can admit stationary solutions under all network conditions; i.e., if the network starts at certain stationary states, then it will stay there all the time. 
These studies on traffic dynamics in a diverge-merge network suggest that the interactions among diverging and merging bottlenecks can lead to unstable traffic patterns.

\bfg\bc
\includegraphics[width=3in]{figure20061015osc.3}
\ec \caption{A diverge-merge network with one O-D pair and two intermediate links: An abstract network}\label{fig:dm2} \efg

In this study, we examine the stability of these stationary states within the framework of kinematic wave models to better understand why the DM network cannot converge to these stationary states when oscillatory solutions occur. 
First we observe that, when one of the two intermediate links is congested and the other not, shock and rarefaction waves travel backward on the congested link and forward on the uncongested, and the two intermediate links form a circular path of information propagation associated with kinematic waves. Based on this observation, we create a section at the downstream boundary of the congested intermediate link and derive a \Po map in the out-flux across the section. Then we demonstrate that the fixed points of the \Po map are related to stationary states in the network. Further we demonstrate that oscillatory traffic dynamics are highly related to unstable fixed points of the \Po map. With the \Po map, we then analyze the stability and bifurcation properties of fixed points of the \Po map with respect to the route choice proportion.

The rest of the paper is organized as follows. 
In Section 2, we present the kinematic wave model of traffic dynamics in the DM network and its stationary solutions.  In Section 3, we derive a \Po map for the network and discuss its finite-time stable fixed points. In Section 4, we study properties of asymptotically stable and unstable fixed points. In Section 5, with examples we demonstrate the bifurcation in the stability of fixed points. In Section 6, we apply the \Po map approach to analyzing traffic patterns in more general road networks. In Section 7, we conclude with discussions and possible future studies.
 
 \section{The kinematic wave model and its stationary solutions} 
For the DM network, we introduce two dummy links at the origin and destination, which are labeled as $r$ and $w$, respectively. Vehicles in this network are categorized into two commodities: vehicles of commodity 1 use link 1, and those of commodity 2 use link 2. On link $a$ ($a\in A=\{r,0,1,2,3,w\}$), a point is denoted by $x_a\in[0,X_a]$, where $X_a$ is the link length. Thus the diverging junction has three coordinates $(0,X_0)\sim (1,0) \sim (2,0)$, and the merging junction has three coordinates $(1,X_1)\sim(2,X_2)\sim (3,0)$. At a point $(a,x_a)$ and time $t$, we define the following quantities of total traffic: density $k_a(x_a,t)$, speed $v_a(x_a,t)$, flow-rate $q_a(x_a,t)$, demand $d_a(x_a,t)$, and supply $s_a(x_a,t)$. In addition, on all links except links 1 and 2, we define $\xi_a(x_a,t)$ as the proportion of commodity 1 vehicles; thus the density of commodity 1 is $\xi_a(x_a,t)k_a(x_a,t)$. Hereafter we omit $(x_a,t)$ from these variables unless necessary.

 \subsection{The kinematic wave model} \label{kwmodel} 
 
\btb
\bc
\begin{tabular}{|l|l|}\hline 
$(a,x_a)$ ($a\in\{r,0,1,2,3,w\}$) & Location $x_a$ on link $a$\\\hline 
$t$  & Time\\\hline 
$X_a$ &Length of link $a$ \\\hline
$l_a\in [0,1]$ & The congested portion of link $a$ \\\hline
$k_a(x_a,t)$ &Density at location $x_a$ and time $t$ on link $a$ \\\hline
$q_a(x_a,t)$ &Flow-rate  \\\hline
$\xi_a(x_a,t)$ &Proportion of commodity 1 vehicles  \\\hline
$\phi_a(x_a,t)$ &Flow-rate (flux) of commodity 1 vehicles \\\hline
$d_a(x_a,t)$ &Demand  \\\hline
$s_a(x_a,t)$ &Supply  \\\hline
$U_a(x_a,t)=(d_a(x_a,t),s_a(x_a,t))$ &Traffic state in the demand-supply space  \\\hline
$C_a$  &Capacity of link $a$ \\\hline
$\beta$ & Merging ratio of link 1\\\hline
$\xi$ & Constant choice proportion of route 1\\\hline
$v_a(t)=q_a(X_a^-,t)$ & Out-flux of link $a$ at time $t$ \\\hline
$u_a(t)=q_a(0^+,t)$ & In-flux of link $a$ at time $t$\\\hline
$\Xi_1$& The set of $\xi$ when information propagates counterclockwise \\\hline
$\Xi_2$& The set of $\xi$ when information propagates clockwise \\\hline
$\tilde \Xi_1$& The set of $\xi$ when the DM network is stationary at SOC-SUC \\\hline
$\tilde \Xi_2$& The set of $\xi$ when the DM network is stationary at SUC-SOC \\\hline
$A_1$ &$\equiv\max\{C_3-(1-\xi) C_0 , C_3-C_2, \beta C_3\}$ defined in \refe{def:A}\\\hline
$A_2$ & $\equiv\max\{C_3-\xi C_0, C_3-C_1, (1-\beta) C_3 \}$ defined in \refe{def:Ap} \\\hline
$A_2'$ & $\equiv C_3-A_2$ defined in \refe{def:A2p}\\\hline
$F\cdot$ & The unified \Po map defined in \refe{pmap} \\\hline
$v^*$ &The fixed point of the \Po map\\\hline
$v^-$ & The smaller periodical point of period 2 of the \Po map\\\hline
$v^+$ & The larger periodical point of period 2 of the \Po map\\\hline
\end{tabular}\caption{A list of notations}\label{listnotations}
\ec
\etb

We assume that all vehicles have the same characteristics, and all links in the DM network are homogeneous. Furthermore, we assume a fundamental diagram at a point $(a,x_a)$ ($a=0,\cdots,3$) and time $t$: $q_a=Q_a(k_a)$, and $v_a=V_a(k_a)=Q_a(k_a)/k_a$ \citep{greenshields1935capacity}. 
Here we assume that $q_a=Q_a(k_a)$ is unimodal and attains its capacity $C_a=Q_a(k_{a,c})$ at the critical density of $k_{a,c}$. 
Then traffic demand and supply are also functions of total density, which are given by
\citep{engquist1980difference,daganzo1995ctm,lebacque1996godunov}
\bsq
\bqn
d_a&=&D_a(k_a)\equiv Q_a(\min\{k_{a,c},k_a\}),\\
s_a&=&S_a(k_a)\equiv Q_a(\max\{k_{a,c},k_a\}).
\eqn
\esq
From the definitions of demand and supply, we can see that the demand-supply pair can uniquely determine density and flow-rate:
\bsq
\bqn
q_a&=&\min\{d_a,s_a\},\\
C_a&=&\max\{d_a,s_a\},\\
k_a&=&R_a(d_a/s_a)\equiv\cas{{ll}D_a^{-1}(C_a d_a/s_a), & d_a\leq s_a\\ S_a^{-1}(C_a s_a/d_a), &s_a\leq d_a}
\eqn
\esq  
Thus we can denote a traffic state at a point by a demand-supply pair $U_a=(d_a,s_a)$.

In the kinematic wave model of traffic dynamics on the DM network, traffic dynamics on link $a$ ($a=0,\cdots,3$) are described by the following LWR model 
\bsq\label{multi-lwr}
\bqn
\pd{k_a}t+\pd{k_a V_a(k_a)} {x_a}&=&0. \label{lwr}
\eqn
To track dynamics of commodity 1 flows on links 0 and 3 we apply a multi-commodity LWR model \citep{lebacque1996godunov}:
\bqn
\pd{\xi_a k_a}t+\pd{\xi_a k_a V_a(k_a)} {x_a}&=&0, \qquad a=0,3. \label{lwr:proportion}
\eqn
\esq  
Then \refe{multi-lwr} constitutes a link-based kinematic wave model of traffic dynamics in the DM network, which is a system of six hyperbolic conservation laws.

In \citep{jin2012network,jin2012statics}, it was shown that the kinematic wave model is well-defined when complemented by the following invariant junction flux functions. Here $\phi_a(x_a,t)$ is the flux of commodity 1 traffic, and $q_a(x_a,t)-\phi_a(x_a,t)$ the flux of commodity 2 traffic.
\ben
\item If a point $(a,x_a)$ has only one upstream link and one downstream link, then the total and commodity 1 fluxes through this point are given by
\bsq \label{link-entropy}
\bqn
q_a(x_a,t)&=&\min\{d_a(x_a^-,t),s_a(x_a^+,t)\},\\
\phi_a(x_a,t) &=&\xi_a(x_a^-,t) q_a(x_a,t), 
\eqn 
\esq
where $x_a^-$ and $x_a^+$ are the upstream and downstream points of $(a,x_a)$, respectively. The flux function can be applied to the origin and destination junctions as follows.
\ben
\item At the origin junction $(0,0)\sim (r,0)$, if the demand at origin, $d_r(0^-,t)$, and the proportion of commodity 1, $\xi_r(0^-,t)$, are given, from \refe{link-entropy} the boundary fluxes are given by 
\bsq
\bqn
q_r(0,t)&=&q_0(0,t)=\min\{d_r(0^-,t),s_0(0^+,t)\},\\
\phi_r(0,t)&=&\xi_r(0^-,t) q_r(0,t).
\eqn
\esq
\item At the destination junction $(3,X_3) \sim (s,0)$, if the supply at the destination, $s_r(0^+,t)$, is given, from \refe{link-entropy} the boundary fluxes are given by 
\bsq
\bqn
q_w(0,t)&=&q_3(X_3,t)=\min\{d_3(X_3^-,t),s_w(0^+,t)\},\\
\phi_w(0,t)&=&\xi_3(X_3^-,t) q_w(0,t).
\eqn
\esq
\een
\item At the diverging junction $(0,X_0)\sim (1,0) \sim (2,0)$, the boundary fluxes are given by
\bsq\label{divergemodel}
\bqn
q_0(X_0,t)&=&\min\{d_0(X_0^-,t),\frac{s_1(0^+,t)}{\xi_0(X_0^-,t)},\frac{s_2(0^+,t)}{1-\xi_0(X_0^-,t)}\},\\
q_1(0,t)&=&\phi_0(X_0,t)=q_0(X_0,t) \xi_0(X_0^-,t),\\
q_2(0,t)&=&q_0(X_0,t) (1-\xi_0(X_0^-,t)).
\eqn
\esq
\item At the merging junction $(1,X_1)\sim (2,X_2)\sim(3,0)$, the boundary fluxes are given by
\bsq \label{mergemodel}
\bqn
q_3(0,t)&=&\min\{d_1(X_1^-,t)+d_2(X_2^-,t),s_3(0^+,t) \},\\
q_1(X_1,t)&=&\phi_3(0,t)=\min\{d_1(X_1^-,t), \max\{s_3(0^+,t)-d_2(X_2^-,t),\beta s_3(0^+,t)\}\},\\
q_2(X_2,t)&=&\min\{d_2(X_2^-,t), \max\{s_3(0^+,t)-d_1(X_1^-,t),(1-\beta) s_3(0^+,t)\},
\eqn
\esq
where $\beta$ is the merging ratio of link 1.
\een

The kinematic wave model, \refe{multi-lwr}, together with the entropy conditions above, defines a semigroup of network hyperbolic conservation laws \citep{bressan1996semigroup}
 in the sense that, given initial conditions in $k_a(x_a,0)$ and $\xi_a(x_a,0)$ and boundary conditions in $d_r(0^-,t)$, $\xi_r(0^-,t)$, and $s_w(0^+,t)$, we can uniquely solve $k_a(x_a,t)$ and $\xi_a(x_a,t)$ at any time. 
Therefore, the network kinematic wave model can be considered as an infinite-dimensional dynamical system.

\subsection{The traffic statics problem and its stationary solutions}

In \citep{jin2012statics}, the traffic statics problem is defined as finding stationary solutions to \refe{multi-lwr} when the origin demand, destination supply, and route choice proportion are all constant: for example, $d_r(0^-,t)=C_0$, $\xi_r(0^-,t)=\xi$, and $s_w(0^+,t)=C_3$. In this case, \refe{multi-lwr} becomes an autonomous, infinite-dimensional system, and traffic dynamics in the road network are determined by the initial conditions on the two intermediate links: $k_1(x_1,0)$ and $k_2(x_2,0)$.

It was found that in stationary states the flow-rates on the two intermediate links are constant: $q_a(x_a,t)=q_a$ ($a=1,2$); and the commodity proportions are constant on links 0 and 3: $\xi_a(x_a,t)=\xi$ ($a=0,3$). Then the flow-rates on links 0 and 3 are the same: $q=q_1+q_2$. 
Since the origin demand equals link 0's capacity, link 0 becomes over-critical (OC) with $k_0(x,t)\geq k_{0,c}$ after a finite period of time; similarly, since the destination supply equals link 3's capacity, link 3 becomes under-critical (UC) with $k_3(x,t)\leq k_{3,c}$ after a finite period of time. Therefore, \refe{multi-lwr} can be simplified into the following two equations:
\bsq \label{simple-lwr}
\bqn
\pd{k_1}t+\pd{k_1 V_1(k_1)} {x_1}&=&0,\\
\pd{k_2}t+\pd{k_2 V_2(k_2)} {x_2}&=&0,
\eqn
\esq 
which are complemented by the following flux functions:
\bsq \label{simple-flux}
\bqn
q_a(x_a,t)&=&\min\{d_a(x_a^-,t),s_a(x_a^+,t)\}, \quad x_a\in(0,X_a), a=1,2\\
q_1(0,t)&=&\min\{\xi C_0 ,s_1(0^+,t),\frac{\xi}{1-\xi} s_2(0^+,t)\} ,\\
q_2(0,t)&=&\min\{(1-\xi) C_0 , \frac{1-\xi}{\xi} s_1(0^+,t), s_2(0^+,t)\},\\
q_1(X_1,t)&=&\min\{d_1(X_1^-,t), \max\{C_3-d_2(X_2^-,t),\beta C_3\}\},\\
q_2(X_2,t)&=&\min\{d_2(X_2^-,t), \max\{C_3-d_1(X_1^-,t),(1-\beta) C_3\}\}.
\eqn
\esq

Then the stationary density on link $a$ can be written as 
\bqn
k_a(x_a,t)&=&(1-I_a(x_a;l_a))R_a(q_a/C_a)+I_a(x_a;l_a)R_a(C_a/q_a), \label{stationarydensity}
\eqn
 where $l_a\in[0,1]$, and the indicator function $I_a(x_a;l_a)=\cas{{ll}0, &x_a\in[0,(1-l_a)X_a)\\1,&x_a\in[(1-l_a)X_a,X_a]}$. On each link, there can be four types of stationary states: (i) when $q_a=C_a$ and any $l_a$, the stationary state is C (critical); (ii) when $q_a<C_a$ and $l_a=0$, the stationary state is SUC (strictly under-critical); (iii) when $q_a<C_a$ and $l_a=1$, the stationary state is SOC (strictly over-critical); and (iv) when $q_a<C_a$ and $l_a\in(0,1)$, the stationary state is ZS (zero shock wave).  

In \citep{jin2012statics} it was shown that stationary solutions to \refe{simple-lwr} with \refe{simple-flux} exist under all network conditions. That is, if the initial condition is a stationary solution, then the network stays at the state. In stationary states, link 0 is always OC at $(C_0,q)$, link 3 is always UC at $(q,C_3)$, and links 1 and 2 can be C, SUC, SOC, or ZS. In \reft{possibless}, solutions of stationary states on links 1 and 2 are listed under all network conditions: from the types of stationary states on both links and the corresponding flow-rates under given network conditions, we can use \refe{stationarydensity} to find the stationary densities on both links. Refer to \citep{jin2012statics} on the derivation of the table.

\btb
\bc
\begin{tabular}{|c|c|c||c|c|}\hline
Capacities& $\xi$ & $\beta$ & Link 1-Link 2 & $q$\\\hline\hline
$C_0<\min\{C_1+C_2, C_3\}$ & $\xi\leq 1-\frac {C_2}{C_0}$& &SUC-C & $C_2/(1-\xi)$\\\hline
&$1-\frac{C_2}{C_0}<\xi<\frac{C_1}{C_0}$& &SUC-SUC & $C_0$\\\hline
&$\xi\geq \frac{C_1}{C_0}$& &C-SUC& $C_1/\xi$ \\\hline\hline
$C_1+C_2\leq \min\{C_0,C_3\}$&$\xi<\frac{C_1}{C_1+C_2}$&&SUC-C& $C_2/(1-\xi)$\\\hline
&$\xi=\frac{C_1}{C_1+C_2}$&&C-C&$C_1/\xi$\\\hline
&$\xi>\frac{C_1}{C_1+C_2}$&&C-SUC&$C_1/\xi$\\\hline\hline
$C_3=C_0<C_1+C_2$&$\xi< 1-\frac {C_2}{C_0}$&&SUC-C& $C_2/(1-\xi)$ \\\hline
&$\xi=1-\frac{C_2}{C_0}$& $\xi<\beta$ &SUC-C& $C_3$\\\hline
&&$\xi\geq\beta$&SUC/SOC/ZS-C&$C_3$\\\hline
&$1-\frac{C_2}{C_0}<\xi<\frac{C_1}{C_0}$&$\xi<\beta$ &SUC-SUC/SOC/ZS& $C_3$\\\hline
&&$\xi=\beta$&SUC/SOC/ZS-SUC/SOC/ZS& $C_3$\\\hline
&&$\xi>\beta$&SUC/SOC/ZS-SUC& $C_3$\\\hline
&$\xi= \frac{C_1}{C_0}$&$\xi\leq\beta$&C-SUC/SOC/ZS&$C_3$\\\hline
&&$\xi>\beta$&C-SUC&$C_3$\\\hline
&$\xi> \frac{C_1}{C_0}$&&C-SUC&$C_1/\xi$\\\hline\hline
$C_3<\min\{C_0,C_1+C_2\}$&$\xi< 1-\frac{C_2}{C_3}$&&SUC-C& $C_2/(1-\xi)$\\\hline
&$\xi= 1-\frac{C_2}{C_3}$&$\xi<\beta$&SUC-C& $C_3$\\\hline
&&$\xi\geq\beta$&SUC/SOC/ZS-C&$C_3$\\\hline
&$1-\frac{C_2}{C_3}<\xi<\frac{C_1}{C_3}$&$\xi<\beta$&SUC-SOC&$C_3$\\\hline
&&$\xi=\beta$&SOC-SUC/SOC/ZS, SUC/ZS-SOC&$C_3$\\\hline
&&$\xi>\beta$&SOC-SUC&$C_3$\\\hline
&$\xi=\frac{C_1}{C_3}$&$\xi\leq\beta$&C-SUC/SOC/ZS&$C_3$\\\hline
&&$\xi>\beta$&C-SUC&$C_3$\\\hline
&$\xi>\frac{C_1}{C_3}$&&C-SUC&$C_1/\xi$\\\hline
\end{tabular}
\caption{Possible stationary states under different network conditions: In the fourth column, the stationary states on the left of a dash are for link 1, and those on the right for link 2. A slash means ``or''. Here $q_1=\xi q$ and $q_2=(1-\xi)q$. \citep{jin2012statics}}\label{possibless}
\ec
\etb

\section{A \Po map of circular information propagation}
When $C_0<\min\{C_1+C_2,C_3\}$ or $C_1+C_2\leq \min\{C_0,C_3\}$; i.e., when the upstream or the middle parts of the DM network are bottlenecks, \reft{possibless} shows that stationary states on links 1 and 2 are both UC. As shown in \citep{jin2009network}, in these cases the stationary states will be reached in a finite time under any initial conditions. Thus, these stationary states are always stable. 

In this study, we focus our discussions on the stability of the kinematic wave model, \refe{simple-lwr} with \refe{simple-flux}, when $C_3\leq C_0$ and $C_3<C_1+C_2$; i.e., when the downstream link 3 imposes a bottleneck for the whole network.
We denote two sets of $\xi$ by 
\bqs
\Xi_1&=&\{\xi: \frac{C_1}{C_3}\leq \xi\leq 1, \m{ or } 1-\frac{C_2}{C_3}<\xi<\frac{C_1}{C_3} \m{ and }\xi\geq \beta; C_3\leq C_0 \m{ and } C_3<C_1+C_2\}\\
\Xi_2&=&\{\xi: 0\leq \xi\leq 1-\frac{C_2}{C_3}, \m{ or } 1-\frac{C_2}{C_3}<\xi<\frac{C_1}{C_3} \m{ and }\xi\leq \beta; C_3\leq C_0 \m{ and } C_3<C_1+C_2\}.
\eqs
The two sets are collectively exhaustive and overlap only when $1-\frac{C_2}{C_3}<\xi<\frac{C_1}{C_3} \m{ and }\xi = \beta$.

\subsection{The derivation of a \Po map}

\begin{figure}\bc
\includegraphics[width=5in]{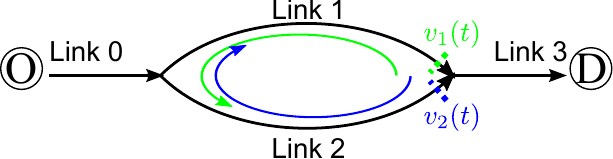}\caption{Circular information propagation in the DM network and two \Po sections}\label{dm2poincare}
\ec
\end{figure}

From \reft{possibless}, when $\xi\in\Xi_1$, links 1 and 2 are stationary at OC and UC, respectively; when $\xi\in\Xi_2$, links 1 and 2 are stationary at UC and OC, respectively. Since shock and rarefaction waves travel upstream in OC states and downstream in UC states, perturbations to a stationary state will circulate with shock or rarefaction waves on the two intermediate links. That is, the two intermediate links are closed by the diverge and the merge to form a circular path for information propagation. In particular, the information propagation direction is counterclockwise when $\xi\in \Xi_1$ and clockwise when $\xi\in \Xi_2$, as shown in \reff{dm2poincare}.
In the following we study the dynamics associated with circular information propagation in terms of \Po maps, which were originally developed to study periodical movements of celestial bodies \citep[][Chapter 10]{wiggins2003chaos}. 
Note that, traditionally, \Po maps were used for ordinary differential equations, but here we derive a \Po map for a system of two partial differential equations in \refe{simple-lwr}. 

First for $\xi\in\Xi_1$, we define a \Po section at the downstream boundary of link 1, $(1,X_1^-)$, as shown in \reff{dm2poincare}.
At $t$, we apply a small perturbation to a stationary state at the point and denote the out-flux by $v_1(t)$. Then the perturbation will travel backward on link 1 to the diverge and then forward on link 2 to the merge. 
When the perturbation reaches the point again after $T$, we obtain $v_1(t+T)$ as a \Po map of $v_1(t)$, $v_1(t+T)=F_1 v_1(t)$, which is a one-dimensional discrete dynamical system. The \Po map can be derived based on the propagation of kinematic waves on the two intermediate links as follows.
\ben
\item  At $t$, a small perturbation is applied on link 1 at $(1,X_1^-)$, so that its out-flux becomes $v_1(t)$. Since link 1 is stationary at OC, we have $d_1(X_1^-,t)=C_1$, and the traffic state at the point becomes $(C_1,v_1(t))$, which propagates upstream. 
\item At $t+T_1$, the perturbed state $(C_1, v_1(t))$ on link 1 reaches the diverge. Since link 2 is stationary at UC, $s_2(0^+,t+T_1)=C_2$. Then from \refe{divergemodel}, the in-flux on link 2 becomes \footnote{Note that here the route choice proportion $\xi$ is always constant.}
  \bqs
  u_2(t+T_1)&=&\min\{(1-\xi) C_0, \frac{1-\xi}{\xi} v_1(t) , C_2\}.
  \eqs
Then the traffic state on link 2 at $(2,0^+)$ becomes $(u_2(t+T_1),C_2)$, which propagates downstream.
  \item At $t=t+T$, \footnote{The exact value of $T$ depends on the length and fundamental diagram of link 1 as well as traffic states, as shown in \citep{jin2009network}.} $(u_2(t+T_1),C_2)$ on link 2 reaches the merge. Then from \refe{mergemodel}, the out-flux on link 1 becomes
  \bqs
  v_1(t+T)&=&\min\{C_1, \max\{C_3-u_2(t+T_1), \beta C_3\}\}.
  \eqs
\een
Therefore, we have
\bqs
v_1(t+T)&=&\min\{C_1, \max\{C_3-\min\{(1-\xi) C_0, \frac{1-\xi}{\xi}  v_1(t) , C_2\}, \beta C_3\}\}\nonumber\\
&=&\min\{C_1, \max\{C_3-(1-\xi) C_0 , C_3-C_2, \beta C_3, C_3-\frac{1-\xi}{\xi}  v_1(t)\}\}\nonumber\\
&=&\min\{C_1,\max\{A_1, C_3-\frac{1-\xi}{\xi}  v_1(t)\}\} \label{discreteds}
\eqs
where 
\bqn
A_1\equiv\max\{C_3-(1-\xi) C_0 , C_3-C_2, \beta C_3\}. \label{def:A}
\eqn
In the derivation, we use some basic properties of min and max operators given in Appendix A. 
Thus, we have the following \Po map 
\bqn
v_1(t+T)&=&F_1 v_1=\min\{C_1, \max\{A_1, C_3-\frac{1-\xi}{\xi} v_1\}\}. \label{def:Fmap}
\eqn

Similarly, for $\xi\in\Xi_2$, we can derive the following \Po map:
\bqn
v_2(t+T')&=&F_2 v_2(t) = \min\{C_2, \max\{A_2,C_3- \frac{\xi}{1-\xi} v_2(t) \}\}, \label{def:Fpmap}
\eqn
where $v_2(t)$ is the perturbed out-flux on link 2 at $(2, X_2^-)$, and 
\bqn
A_2\equiv\max\{C_3-\xi C_0, C_3-C_1, (1-\beta) C_3 \}. \label{def:Ap}
\eqn
Comparing \refe{def:Fmap} and \refe{def:Fpmap}, we can see that they are symmetric in the sense that they are equivalent if $v_1$, $\xi$, $C_1$, and $\beta$ are swapped with $v_2$, $1-\xi$, $C_2$, and $1-\beta$ respectively. 

Note that a discrete dynamical system similar to \refe{def:Fmap} was first derived in \citep{jin2009network}, but was not thoroughly analyzed as \Po maps. Here the new \Po maps, \refe{def:Fmap} and \refe{def:Fpmap}, are more complete, since they apply to a more general merging model \refe{mergemodel} and more network conditions.

If we define a new variable $v(t)=\cas{{ll}v_1(t),&\xi\in \Xi_1\\C_3-v_2(t),&\xi \in \Xi_2}$, then we can obtain a unified \Po map: $v(t+T)=F v(t)$, where
\bqn
F v&\equiv&\cas{{ll}F_1 v, &\xi\in\Xi_1\\C_3-F_2(C_3-v),&\xi\in\Xi_2}\nonumber\\
&=&\cas{{ll}\min\{C_1, \max\{A_1, C_3-\frac{1-\xi}{\xi} v\}\},&\xi\in \Xi_1\\
\max\{C_3-C_2,\min\{A_2',\frac{\xi}{1-\xi}(C_3-v)\}\}, &\xi\in\Xi_2 \label{pmap}
}
\eqn  
where 
\bqn
A_2'\equiv C_3-A_2=\min\{\xi C_0,C_1,\beta C_3\}. \label{def:A2p}
\eqn

From \reft{possibless}, it can be easily verified that the fixed point of \refe{pmap} has the same flow-rate as in the corresponding stationary states. That is, we have the following theorem.
\begin{theorem} \label{thm:fixedpoint}
When $\xi\in \Xi_1$, the fixed point of \refe{pmap} is $v^*=\xi q$, where $q$ is given in \reft{possibless} under the same network conditions. 
Namely, the fixed point is
\bqn
v^*&=&\cas{{ll} C_1, & \frac{C_1}{C_3} \leq \xi \leq 1;\\ \xi C_3, & \xi\in(1-\frac{C_2}{C_3},\frac{C_1}{C_3}) \m{ and } \xi \geq \beta.} \label{Fmap:fp}
\eqn
 When $\xi \in \Xi_2$, the fixed point of \refe{pmap} is $v^*=C_3-(1-\xi) q$, where $q$ is given in \reft{possibless} under the same network conditions.
 Namely, the fixed point is
 \bqn
 v^*&=&\cas{{ll} C_3-C_2, & 0\leq \xi \leq 1-\frac{C_2}{C_3};\\ \xi C_3, & \xi\in(1-\frac{C_2}{C_3},\frac{C_1}{C_3})\m{ and } \xi\leq \beta.} \label{Fpmap:fp}
 \eqn 
\end{theorem}

\subsection{Finite-time stable fixed points of the \Po  map}
For a discrete dynamical system, \refe{pmap}, a fixed point is finite-time stable, if and only if the dynamical system, starting from any initial conditions, converges to the fixed point in a finite number of time-steps \citep{haimo1986finite,bhat2000finite}.
Then we have the following theorem, whose proof is given in Appendix B.
\begin{theorem} \label{thm:fts} For DM networks with $\xi\in \Xi_1$, \refe{pmap} is finite-time stable at $v^*=C_1$ when  $\frac {C_1}{C_3}\leq \xi\leq 1$; or at $v^*=\xi C_3$ when $\xi\in(1-\frac{C_2}{C_3},\frac{C_1}{C_3})$ and $\xi= \beta$; or at $v^*=\xi C_3$ when $C_3=C_0$, $\xi\in(1-\frac{C_2}{C_3},\frac{C_1}{C_3})$, and $\xi>\beta$. Similarly, for $\xi\in \Xi_2$, \refe{pmap} is finite-time stable at $v^*=C_3-C_2$ when $0\leq\xi \leq 1-\frac{C_2}{C_3}$; or at $v^*=\xi C_3$ when $C_3=C_0$, $\xi\in(1-\frac{C_2}{C_3},\frac{C_1}{C_3})$, and $\xi< \beta$.
\end{theorem}

In \citep{jin2009network}, it was shown that the kinematic wave model \refe{simple-lwr} also reaches stationary states in a finite period of time under the corresponding network conditions. This demonstrates the consistency between the \Po map \refe{pmap} and the original kinematic wave model in finite-time stable traffic dynamics.

\section{Stability of SOC-SUC and SUC-SOC stationary states}
From Theorem \ref{thm:fts}, the DM network converge to stationary states in a finite period of time except for $C_3<\min\{C_0, C_1+C_2\}$ and $\xi\in \tilde \Xi_1$ or $\xi\in \tilde \Xi_2$, where
\bqs
\tilde \Xi_1&=&\{\xi\in(1-\frac{C_2}{C_3}, \frac{C_1}{C_3}) \m{ and } \xi>\beta; C_3<\min\{C_0, C_1+C_2\}\}\\
\tilde \Xi_2&=&\{\xi\in(1-\frac{C_2}{C_3},\frac{C_1}{C_3}) \m{ and } \xi<\beta; C_3<\min\{C_0, C_1+C_2\}\}
\eqs
From \reft{possibless}, we can see that the network is stationary at SOC-SUC and SUC-SOC for $\xi\in \tilde \Xi_1$ and $\tilde \Xi_2$, respectively. In both cases, the fixed point of \refe{pmap} is $v^*=\xi C_3$. 
In this section we focus on the stability of the fixed point.

\subsection{Local stability analysis with Lyapunov's first method}
When $\xi\in \tilde \Xi_1$, we denote a small perturbation of $v$ around the fixed point by $\tilde v=v-\xi C_3$. Since $v\approx \xi C_3$, we have $C_3-\frac{1-\xi}{\xi} v\approx \xi C_3$, which is greater than $A_1$ and smaller than $C_1$. Thus \refe{pmap} is equivalent to a local linear map:
\bqn
\tilde F \tilde v&=&-\frac{1-\xi}{\xi} \tilde v, \label{stabilityds}
\eqn 
which leads to $\tilde F^n \tilde v=(-\frac{1-\xi}{\xi})^n \tilde v$.
Since $\frac{1-\xi}{\xi}>0$ for $\xi\in \tilde \Xi_1$, the map is always oscillatory.  The local stability of \refe{pmap} depends on the magnitude of $\frac{1-\xi}{\xi}$. Following the first Lyapunov method  \citep{lasalle1976stability,slotine1991applied,galor2004introduction}, we obtain the following stability properties of the original \Po map, \refe{pmap}.  
\begin{theorem}\label{lyapunov1}
For a DM network with $\xi\in \tilde \Xi_1$, which admits a SOC-SUC stationary state, the stability of the fixed point of the \Po map, \refe{pmap}, is as follows:
(i) When $\frac{1-\xi}{\xi}\geq 1$; i.e., when $\xi\leq \frac 12$,  \refe{stabilityds} is unstable, and, therefore, \refe{pmap} is locally unstable. In this case, the SOC-SUC stationary state is unstable, and PPO (persistent periodic oscillatory) traffic patterns develop in the network; 
(ii) When $0\leq \frac{1-\xi}{\xi}<1$; i.e., when $\xi> \frac 12$, \refe{stabilityds} is asymptotically stable, and, therefore,  \refe{pmap} is asymptotically stable. In this case, the SOC-SUC stationary state is asymptotically stable, and  DPO (damped periodic oscillatory) traffic patterns develop in the network.
\end{theorem}

Similarly, we can obtain the following corollary for fixed points corresponding to SUC-SOC stationary states when $\xi\in\tilde \Xi_2$.
\begin{corollary}\label{lyapunov1cor}
For a DM network with $\xi\in \tilde \Xi_2$, which admits a SUC-SOC stationary state, the stability of the fixed point of the \Po map, \refe{pmap}, is as follows:
(i) When $\frac{\xi}{1-\xi}\geq 1$; i.e., when $\xi\geq \frac 12$,  \refe{pmap} is locally unstable. In this case, the SUC-SOC stationary state is unstable, and the DM network has PPO solutions;
(ii) When $0\leq \frac{\xi}{1-\xi}<1$; i.e., when $\xi< \frac 12$, \refe{pmap} is asymptotically stable. In this case, the SUC-SOC stationary state is asymptotically stable, and the DM network has DPO solutions.
\end{corollary}

Here we obtain the necessary and sufficient conditions for the occurrence of both PPO and DPO solutions in a DM network, which are consistent with those  in \citep{jin2009network}. 

\subsection{Periodical points of the \Po map}
Furthermore, as shown in \citep{jin2009network}, the PPO solutions tend to have constant oscillation magnitudes in both densities and flow-rates. In this  subsection, we will relate these magnitudes to periodical points of the \Po map, \refe{pmap} with $\xi\in \tilde \Xi_1$ or $\tilde \Xi_2$. 
For a discrete \Po  map \refe{pmap}, if $F^n v=v$ and $F^k v\neq v$ for $k<n$, then $v$ is a periodical point of period $n$ \citep[][Chapter 4]{holmgren1996first}.  In this sense, the fixed point of the \Po map, $v^*=\xi C_3$, is a periodical point of period $1$. 

When $\xi\in\tilde \Xi_1$, we first attempt to find periodical points of period 2, which satisfy $F^2 v=v$; i.e., 
\bqn
v&=&F^2 v\equiv\min\{C_1, \max\{A_1, C_3-\frac{1-\xi}{\xi} \min\{C_1, \max\{A_1, C_3-\frac{1-\xi}{\xi} v\}\}\}\}. \label{eqn:period2}
\eqn
Obviously, the fixed point $v^*=\xi C_3$ is a solution of \refe{eqn:period2}.

When the fixed point of \refe{pmap} is asymptotically stable; i.e., when $\xi>\frac 12$, we have the following lemma, whose proof is given in Appendix C.
\begin{lemma} \label{lemma:p1}
For a DM network with $\xi\in\tilde \Xi_1$ and $\xi>\frac 12$, \refe{eqn:period2} has no solutions different from $v^*=\xi C_3$; i.e., the \Po  map, \refe{pmap}, has no periodical points of period 2 when \refe{pmap} is asymptotically stable or when the network has DPO solutions.
\end{lemma}

When the fixed point of \refe{pmap} is unstable; i.e., when $\xi\leq \frac 12$, we have the following lemma, whose proof is also given in Appendix C.
\begin{lemma}\label{lemma:p2}
For a DM network with $\xi\in\tilde \Xi_1$ and $\xi\leq \frac 12$, the \Po  map, \refe{pmap}, has two periodical points of period 2, $v^-$ and $v^+$:
\bsq
\bqn
v^-&=&\max\{A_1,C_3-\frac{1-\xi}{\xi} C_1\}<v^*,\\
v^+&=&\min\{C_1, C_3-\frac{1-\xi}{\xi} A_1\}>v^*.
\eqn
\esq
That is, $v^-$ and $v^+$ are two solutions of \refe{eqn:period2}, which are different from the fixed point. 
In addition,
\bsq
\bqn
F v^-&=&v^+,\\
F v^+&=&v^-,
\eqn
\esq
Thus $\{v^*, v^-,v^+\}$ is an invariant set of the \Po  map \citep{lasalle1976stability}.
When $\xi=\frac 12$, any $v\in[v^-,v^*)\cup (v^*,v^+]$ is a periodical point of period 2.
\end{lemma}

Furthermore, the following lemma, whose proof is also given in Appendix C, shows that there are no periodical points of period 4 for \refe{pmap} when $\xi\leq \frac 12$.
\begin{lemma} \label{lemma:p3} For a DM network with $\xi\in\tilde\Xi_1$ and $\xi\leq\frac 12$, there are no periodical points of period 4 for the \Po  map \refe{pmap}.
\end{lemma}

From Lemmas \ref{lemma:p1}, \ref{lemma:p2}, and \ref{lemma:p3}, we then have the following theorem.
\begin{theorem} \label{thm:period}
For a DM network with $\xi\in\tilde\Xi_1$, which has a SOC-SUC stationary state, the \Po  map, \refe{pmap}, has one periodical point of period 1, $v^*=\xi C_3$. When $\xi < \frac 12$, it has two periodical points of period 2, $v^-=\max\{A_1,C_3-\frac{1-\xi}{\xi} C_1\}<v^*$ and $v^+=\min\{C_1, C_3-\frac{1-\xi}{\xi} A_1\}>v^*$; when $\xi=\frac 12$, it has an infinite number of periodical points of period 2 between $v^-=\max\{A_1,C_3- C_1\}$ and $v^+=\min\{C_1, C_3- A_1\}$. The map has no other periodical points.
\end{theorem}
{\em Proof}. According to Sarkovskii's theorem \citep[][Chapter 5]{holmgren1996first}, if a discrete map has no periodical points of period 2, it only has a fixed point of period 1. Then from Lemma \ref{lemma:p1}, the \Po  map \refe{def:Fmap} has only one fixed point of period 1 when $\xi>\frac 12$.

Sarkovskii's theorem also states that, if a discrete map has no periodical points of period 4, then it only has periodical points of periods 1 and 2. Then from Lemma \refe{lemma:p2}, the \Po  map \refe{def:Fmap} has only one fixed point of period 1, two periodical points of period 2 when $\xi< \frac 12$, and an infinite number of periodical points of period 2 when $\xi=\frac 12$. \eop

Similar to Theorem \ref{thm:period}, the following corollary states periodical points of the \Po map when $\xi\in\tilde \Xi_2$.
\begin{corollary} \label{Gperiod}
For a DM network with  $\xi\in\tilde \Xi_2$, which has a SUC-SOC stationary state, the \Po  map, \refe{pmap}, has one periodical point of period 1, $v^*=\xi C_3$. 
When $\xi>\frac 12$, it has two periodical points of period 2, $v^-=\max\{C_3-C_2,\frac{\xi}{1-\xi} A_2'\}<v^*$ and $v^+=\min\{A_2',\frac{\xi}{1-\xi}C_2\}>v^*$; when $\xi=\frac 12$, any point $v\in[v^-,v^*)\cup (v^*, v^+]$ is a periodical point of period 2. The map has no other periodical points.
\end{corollary}

In addition, we have the following theorem regarding chaotic solutions to the \Po map, \refe{pmap}.
\begin{corollary}
For a DM network with $\xi\in\tilde \Xi_1$ or $\xi\in\tilde \Xi_2$, the \Po  map, \refe{pmap}, has no chaotic solutions.
\end{corollary}
{\em Proof}. According to \citep{li1975period,burns2011sharkovsky}, chaotic solutions exist if and only if when there is a periodical point of period 3. Thus from Theorem \ref{thm:period} and Corollary \ref{Gperiod}, there exist no chaotic solutions for the \Po map. \eop

In the following we consider two examples. 
First, we consider the DM network studied in \citep{jin2003dissertation,jin2005paramics,jin2009network}: $C_0=3$, $C_1=1$, $C_2=2$, $C_3=2$, and $\beta=\frac{C_1}{C_1+C_2}=\frac 13$. In this network, $C_3<\min\{C_0,C_1+C_2\}$,  $A_1=\max\{3\xi-1,\frac 23\}$, and $A_2'=\min\{3\xi,\frac 23\}$. The \Po map \refe{pmap} can be written as
\bqs
F v&=&\cas{{ll} \min\{1,\max\{3\xi-1,\frac 23,2-\frac{1-\xi}\xi v\}\}, & \frac 13\leq \xi\leq 1\\\max\{0,\min\{3\xi,\frac23,\frac{\xi}{1-\xi}(2-v)\}\},&0\leq\xi\leq\frac 13}
\eqs
\ben
\item From Theorem \ref{thm:fts}, the \Po map is finite-time stable at $v^*=C_1=1$ when $\frac 12\leq\xi\leq 1$; or at $v^*=\xi C_3=\frac 23$ when $\xi=\frac 13$; or at $v^*=C_3-C_2=0$ when $\xi=0$.  
\item From Theorem \ref{lyapunov1} and Corollary \ref{lyapunov1cor}, when $\xi\in(\frac 13, \frac 12)$, the map is locally unstable at $v^*=\xi C_3=2\xi$, which corresponds to a SOC-SUC stationary state, since $\xi\leq \frac 12$; when $\xi\in(0, \frac 13)$, the map is asymptotically stable at $v^*=\xi C_3=2\xi$, which corresponds to a SUC-SOC stationary state. 
Comparing with results in \citep{jin2009network}, we can verify that PPO traffic patterns arise in an initially empty DM network when the stationary states are unstable, and DPO traffic patterns arise when the stationary states are asymptotically stable. Therefore, unstable stationary states can never be reached in an initially empty network, and asymptotically stable stationary states can be reached only after a long time. Note that this conclusion does not contradict that in \citep{jin2012statics}, where it was shown that stationary states always exist for the DM network. 
\item When $\xi\in(\frac 13,\frac 12)$, $A_1=\max\{2-3(1-\xi),0,\frac 23\}=\max\{3\xi-1,0,\frac 23\}=\frac 23$. Thus from Theorem \ref{thm:period}, there are two periodical points of period 2 for the \Po map:
\bqs
v^-&=&\max\{A,2-\frac{1-\xi}\xi \}=\max\{\frac 23,3-\frac 1\xi\},\\
v^+&=&\min\{1,2-\frac{1-\xi}\xi A\}=\min\{1,\frac 83-\frac 23 \frac 1\xi\}.
\eqs
For example, when $\xi=0.45$, $v^-=\frac 79$, and $v^+=1$. These two values correspond to the lower and upper magnitudes of flow-rates on link 1 in the PPO solutions as observed in \citep[][Section 7.3]{jin2003dissertation} and \citep{jin2009network}. Therefore, the periodical points of periodic 2 determine the oscillation magnitudes of PPO solutions in the DM network.

\een

In the second example, we consider a DM network with $C_3<\min\{C_0,C_1+C_2\}$, $C_1=C_2$, and the fair merging rule with $\beta=C_1/(C_1+C_2)=\frac 12$. 
From Theorem \ref{thm:fts}, we can find all finite-time stable fixed points of the \Po map, \refe{pmap}.
This network has SUC-SOC stationary states when $\xi\in(1-\frac{C_2}{C_3},\frac{C_1}{C_3})$ and $\xi<\frac 12$, and SOC-SUC stationary states when $\xi\in(1-\frac{C_2}{C_3},\frac{C_1}{C_3})$ and $\xi>\frac 12$.
  But from Theorem \ref{lyapunov1} and Corollary \ref{lyapunov1cor}, the corresponding stationary state are always asymptotically stable. That is, the DM network has DPO solutions, but no PPO solutions. This is again consistent with numerical simulations in \citep[][Section 7.2]{jin2003dissertation}.

\section{Bifurcation of stationary states with respect to route choice proportions}

In this section, we consider the impacts of route choice proportions on traffic dynamics in the following DM network: $C_0=3$, $C_1=1.5$, $C_2=2$, $C_3=2.5$, and $\beta=0.3$. In this network, $C_3<\min\{C_0,C_1+C_2\}$, $A_1=\max\{3 \xi -0.5, 0.75\}$, and $A_2'=\min\{3\xi, 0.75\}$. The \Po map, \refe{pmap}, can be written as
\bqs
Fv&=&\cas{{ll}\min\{1.5, \max\{3 \xi -0.5, 0.75, 2.5-\frac{1-\xi}{\xi}v\}\}, & \xi \in[0.3,1]\\\max\{0.5, \min\{3\xi, 0.75, \frac{\xi}{1-\xi} (2.5-v)\}\},& \xi\in[0,0.3]}
\eqs

Then the stability of the \Po map is as follows:
\ben
\item From Theorem \ref{thm:fts}, the \Po map is finite-time stable at $v^*=C_1=1.5$ when $\xi\in[0.6,1]$; or at $v^*=\beta C_3=0.75$ when $\xi=0.3$; or $v^*=C_3-C_2=0.5$ when $\xi\in[0,0.2]$. 
\item From Theorem \ref{lyapunov1} and Corollary \ref{lyapunov1cor}, when $\xi\in(0.3, 0.5]$, the map is locally unstable at $v^*=\xi C_3=2.5\xi$, which corresponds to a SOC-SUC stationary state; when $\xi\in(0.5,0.6)$, it is asymptotically stable at $v^*=\xi C_3=2.5\xi$, which also corresponds to a SOC-SUC stationary state; when $\xi\in(0.2, 0.3)$, the map is asymptotically stable at $v^*=\xi C_3=2.5\xi$, which corresponds to a SUC-SOC stationary state. 
\item From Theorem \ref{thm:period}, there are two periodical points of period 2 for the map when $\xi\in(0.3,0.5)$:
\bqs
v^-&=&\max\{3\xi-0.5, 0.75, 2.5-\frac{1-\xi}{\xi} 1.5 \}\\
v^+&=&\min\{1.5, 2.5-\frac{1-\xi}{\xi} 0.75, 3\xi-0.5+0.5\frac{1-\xi}{\xi}\}.
\eqs
When $\xi=0.5$, any point $v\in[v^-,v*)\cup(v^*,v^+]=[1,1.25)\cup(1.25,1.5]$ is a periodical point of period 2. 
\een 

In the following we consider two numerical examples:
\ben
\item When $\xi=0.55$ with initial $v=1.1$, solutions of $F^n v$ and the orbit of the \Po  map \refe{pmap} are shown in \reff{Fmap1}. 
Since $\xi>0.5$, the \Po map is asymptotically stable, and $F^n v$ converges to $v^*=\xi C_3=1.375$.
\item When $\xi=0.4$, $A_1=0.75$, and solutions of $F^n v$ and the orbit of the \Po map \refe{pmap} are shown in \reff{Fmap2}.
Since $\xi<0.5$, the map is locally unstable and oscillates between two periodical points: $v^-=0.75$ and $v^+=1.375$.
\een
We can clearly see that the numerical results are consistent with the theoretical predictions.
 
\bfg\bc $\ba{c@{\hspace{0.3in}}c}
\includegraphics[width=3.2in]{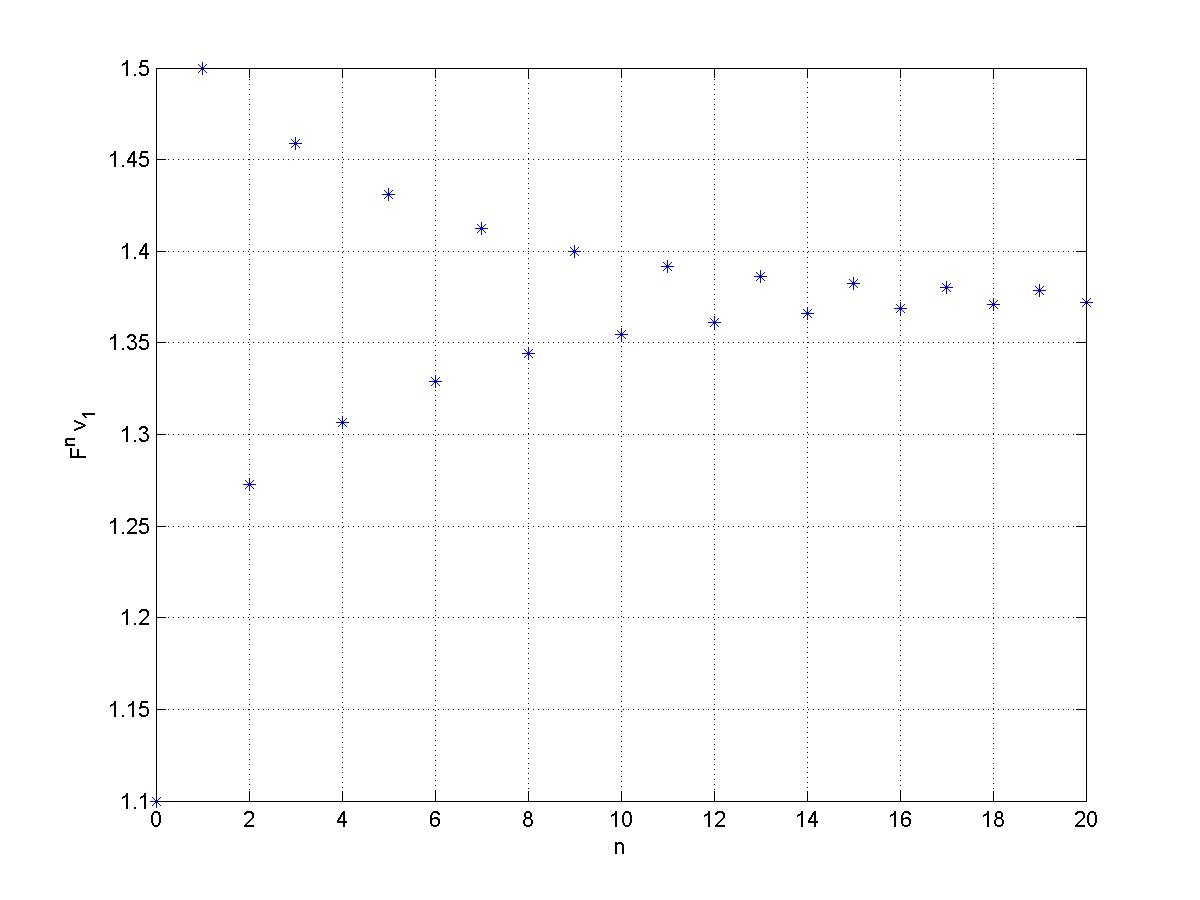} &
\includegraphics[width=3.2in]{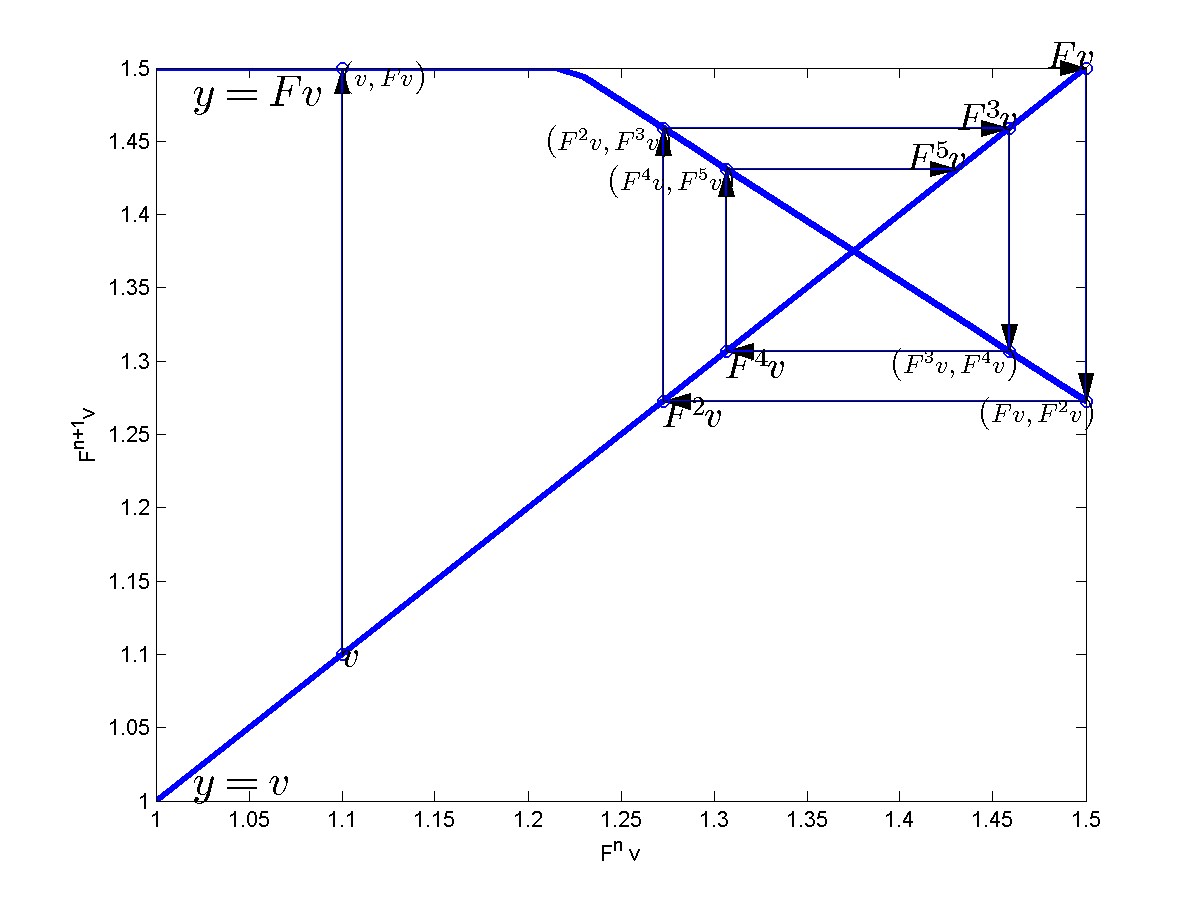} \\
\multicolumn{1}{c}{\mbox{\bf (a)}} &
    \multicolumn{1}{c}{\mbox{\bf (b)}}
\ea$ \ec \caption{Solutions of $F^n v$ and the orbit of the \Po  map \refe{def:Fmap} for $C_0=3$, $C_1=1.5$, $C_2=2$, $C_3=2.5$, $\beta=0.3$, and  $\xi=0.55$ with initial $v=1.1$}\label{Fmap1} \efg

\bfg\bc $\ba{c@{\hspace{0.3in}}c}
\includegraphics[width=3.2in]{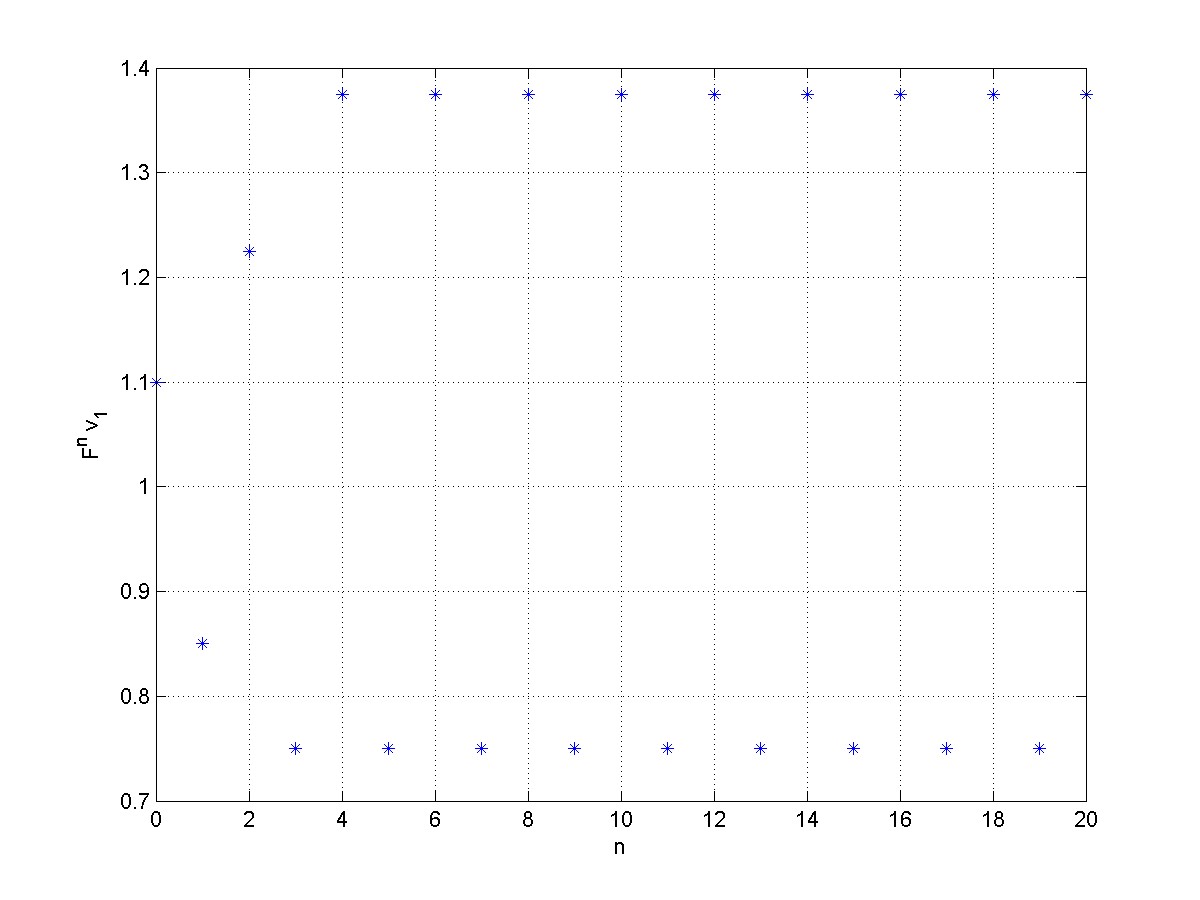} &
\includegraphics[width=3.2in]{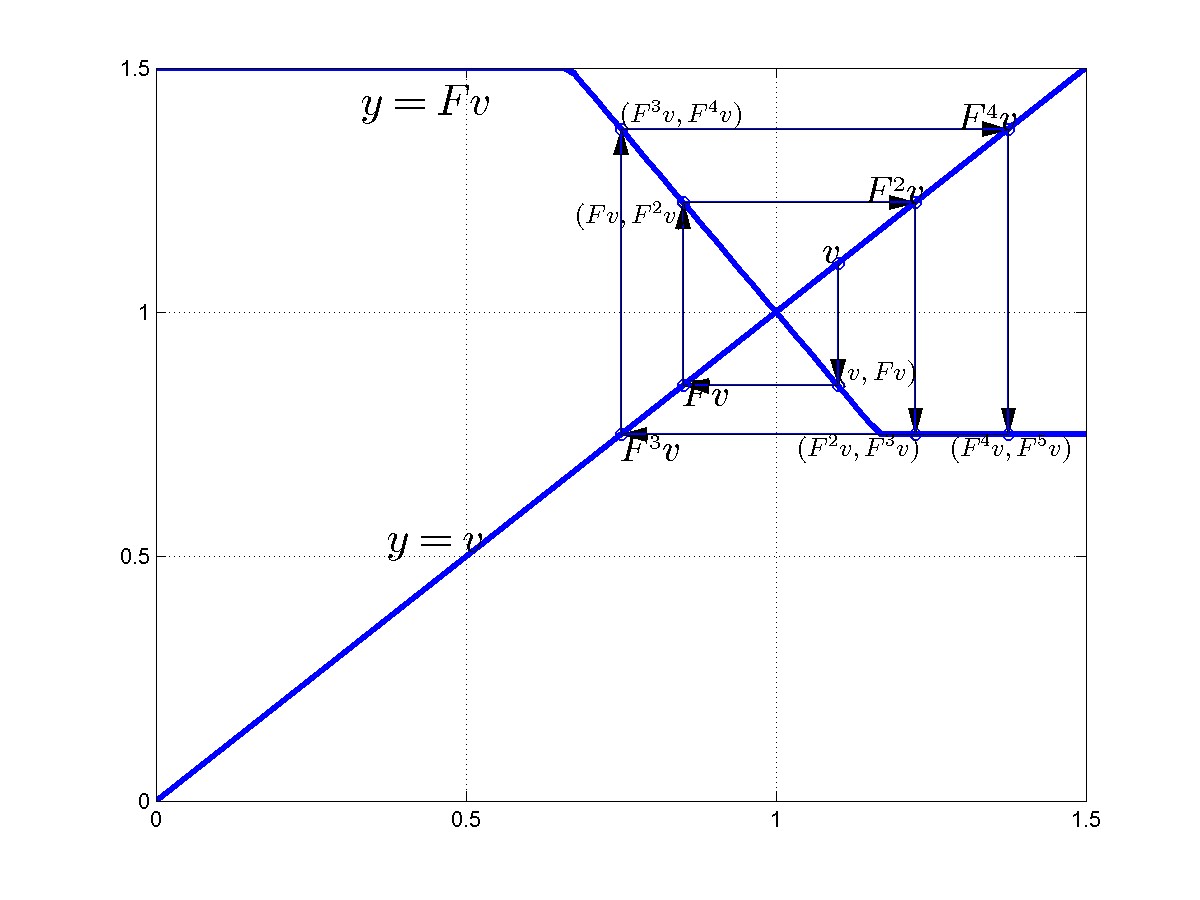} \\
\multicolumn{1}{c}{\mbox{\bf (a)}} &
    \multicolumn{1}{c}{\mbox{\bf (b)}}
\ea$ \ec \caption{Solutions of $F^n v$ and the orbit of the \Po  map \refe{def:Fmap} for $C_0=3$, $C_1=1.5$, $C_2=2$, $C_3=2.5$, $\beta=0.3$, and  $\xi=0.4$ with initial $v=1.1$}\label{Fmap2} \efg

We can see that the \Po map \refe{pmap} always has one fixed point, but its stability property of can change with the route choice proportion, $\xi$, in a DM network when $C_3<\min\{C_0, C_1+C_2\}$: it is finite-time stable for $\xi\in[0,0.2]\cup \{0.3\}\cup[0.6,1]$, asymptotically stable for $\xi\in(0.2,0.3)\cup(0.5,0.6)$, and unstable for $\xi\in(0.3,0.5]$. 
Therefore, the Hopf bifurcation phenomenon occurs with respect to the route choice proportion $\xi$, which is the bifurcation parameter \citep[][Chapter 7]{holmgren1996first}. Here a pair of periodical points of period two form a one-dimensional limit cycle. 
The bifurcation diagram of \refe{pmap} is shown in \reff{bifurcation1}. In the figure, all fixed points are shown on the solid and thick piecewise linear curves, finite-time stable fixed points are marked by circles, and the points on the dashed lines are for periodical points of period 2. In addition, the arrows show the direction of iterations. Thus we can see that fixed points for $\xi\in(0.2,0.3)$ and $\xi\in(0.5,0.6)$ are stable, those for  $\xi\in(0.3,0.5]$ are unstable, but the periodical points are stable.

Since the route choice proportions are determined by drivers' route choice behaviors in a road network, the bifurcation study suggests that a traffic system's stability can be significantly impacted by drivers' choices of routes.

\bfg\bc
\includegraphics[width=6in]{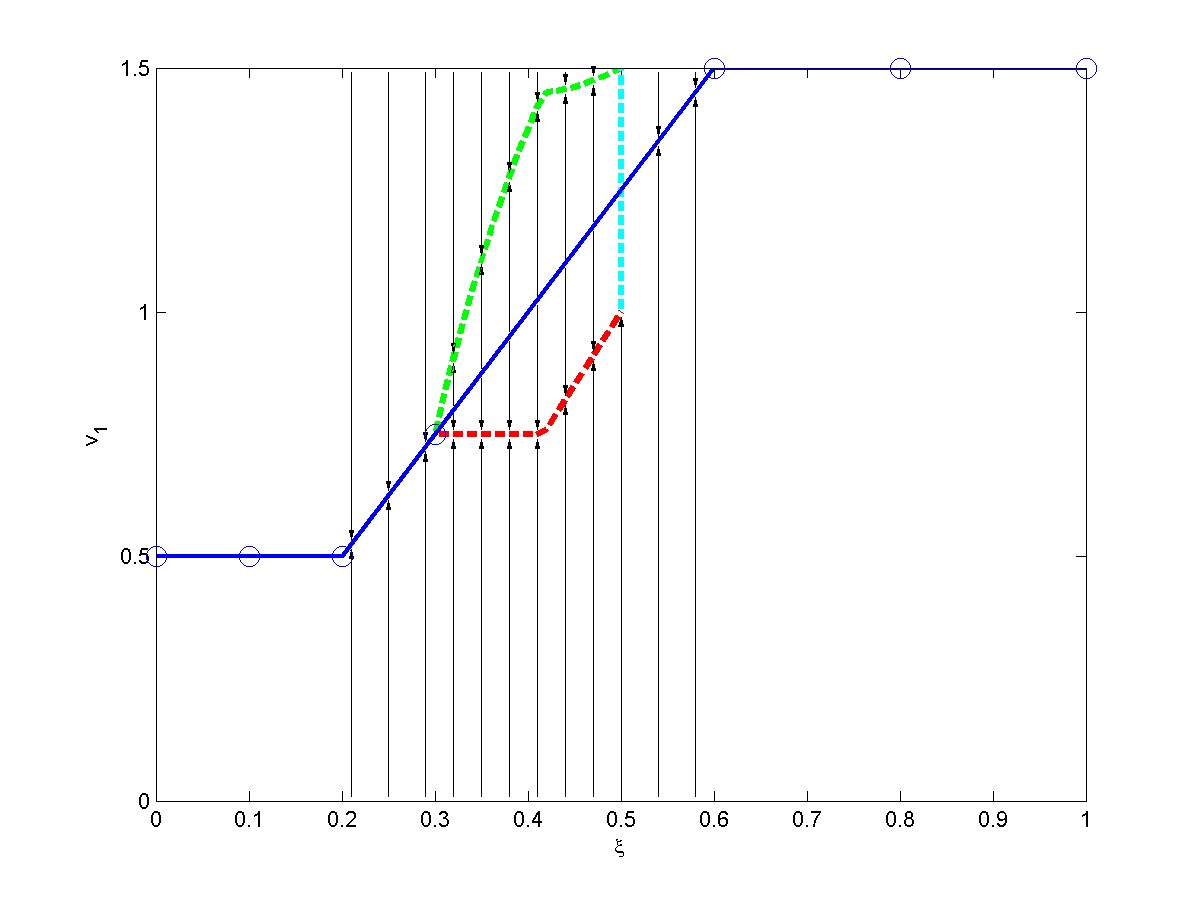}\caption{The bifurcation diagram of the \Po  map \refe{pmap}: $C_0=3$, $C_1=1.5$, $C_2=2$, $C_3=2.5$, and $\beta=0.3$. Here $\xi$ is the bifurcation parameter.}\label{bifurcation1}
\ec\efg

\commentout{

Periodicity follows from Excercise 4.4 of Chapter 1 \citep{lasalle1976stability}.

Three types of stability: zeroth-order stable, asymptotically stable, limit cycle

\del{to add simulation results here?}

\citep{kuznetsov2004elements}: Definition 2.11 ``The apperance of a topologically nonequivalent phase portrait under variation of parameters is called a bifurcation.'' ``bifurcation parameter value'', ``local bifurcations as bifurcations of equilibria or fixed points''
flip bifurcation? section 4.4

\add{This means that the mapping from initial conditions to final solutions is not continuous?}

\add{strong stability?}

}

\section{Application of the \Po map approach to analyzing traffic dynamics in general network structures}
From the analyses in the preceding sections, an insight  is that the \Po map approach can be used to analyze traffic dynamics when shock and rarefaction waves propagate in a circular pattern, which is caused by the bottleneck effects of diverging and merging junctions. Another insight is that persistent periodic oscillations can occur when the \Po map is unstable locally. In this section, we further apply the \Po map to analyzing traffic dynamics in more general networks.

\subsection{$(DM)^n$ networks}
\begin{figure} \bc
\includegraphics[height=2.5in]{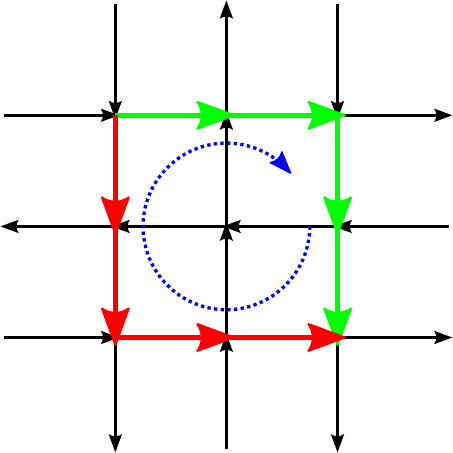} 
\caption{A $3\times 3$ grid network that embeds the DM network} \label{3by3dm2} \ec 
\end{figure}

First, we expect unstable traffic patterns can be observed in networks that embed a DM network. An example is the $3\times 3$ grid network of one-way roads shown in \reff{3by3dm2}. It is possible that, under certain demand patterns, turning proportions, and traffic signals, traffic dynamics are dominated by those in the DM network. If green and red links carry free and congested flow respectively, there exists a circular information path as shown by the dashed line, whose arrows represent the information propagation direction. Obviously we can apply the same \Po map, \refe{pmap}, to demonstrate that traffic can become unstable in the network. In this sense, the DM network structure is a sufficient condition for the existence of unstable traffic dynamics in a road network.

\begin{figure} \bc
$\ba{c@{\hspace{0.3in}}c}
\includegraphics[height=2.5in]{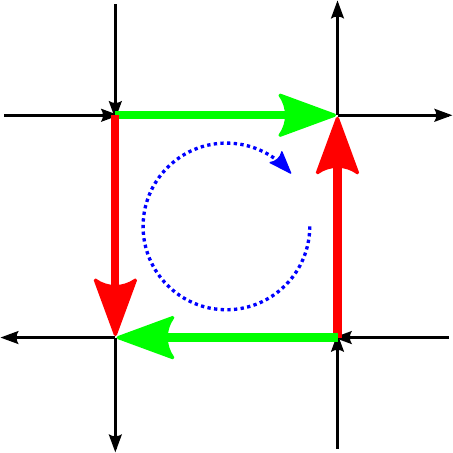} &
\includegraphics[height=2.5in]{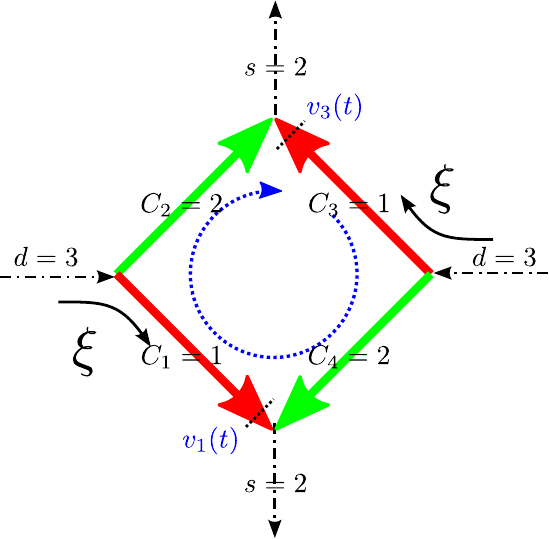} \\
\multicolumn{1}{c}{\mbox{\bf (a)}} &
    \multicolumn{1}{c}{\mbox{\bf (b)}}
\ea$
\caption{A $2\times 2$ grid network (a) and the $(DM)^2$ network (b)} \label{2by2dm2squared} \ec 
\end{figure}

Second, we apply the \Po map approach to analyzing traffic dynamics in a $2\times 2$ grid network, shown in \reff{2by2dm2squared}(a), in which information propagates in a circular fashion when red links are congested and green links not. We assume that traffic dynamics are dominated by the subnetwork shown in \reff{2by2dm2squared}(a), which has two diverge junctions and two merge junctions. We refer to the network in \reff{2by2dm2squared}(b) as a $(DM)^2$ network. Here we consider a symmetric case shown in  in \reff{2by2dm2squared}(b), in which all links have the same length, free-flow speed, and shock wave speed in congested traffic, the origin demands $d=3$, destination supplies $s=2$, the capacities of four links are 1, 2, 1, and 2, respectively, and the turning proportion is $\xi$ to links 1 and 3. We denote the out-fluxes of links 1 and 3 by $v_1(t)$ and $v_3(t)$, respectively. Since links 1 and 3 are congested, and links 2 and 4 not for $\frac 13\xi\frac 12$, we can obtain the following \Po map:
\bsq \label{dm2squared_pm}
\bqn
v_1(t+T)&=&\min\{1,2-\frac{1-\xi}\xi v_3(t)\},\\
v_3(t+T)&=&\min\{1,2-\frac{1-\xi}\xi v_1(t)\},
\eqn
\esq
where $T$ is determined by the link length, free-flow speed, and shock wave speed in congested traffic. If initial conditions are also symmetric, then $v_1(t)=v_3(t)$, and \refe{dm2squared_pm} is equivalent to the \Po map for the DM network:
\bqs
v_1(t+T)&=&\min\{1,2-\frac{1-\xi}\xi v_3(t)\},
\eqs
which is a special case of \refe{pmap} discussed in Section 4.2. In this case, the network has one stationary state $v_1^*=v_3^*= 2\xi$, which is unstable and oscillates between $v^-=\max\{\frac 23, 3-\frac 1\xi\}$ and $v^+=\max\{1,\frac 83-\frac 23 \frac 1\xi\}$. However, if the initial conditions are asymmetric, \refe{dm2squared_pm} becomes
\bqs
v_1(t+2T)&=&\min\{1,2-\frac{1-\xi}\xi \min\{1,2-\frac{1-\xi}\xi v_1(t)\}\},
\eqs
for which $v_1^*=2\xi$ is still a fixed point (stationary state). If we apply a perturbation $\tilde v_1(t)$ to $v_1^*$, we have
\bqs
\tilde v_1(t+2T)&=&(\frac{1-\xi}\xi)^2 \tilde v_1(t),
\eqs
which is unstable, since $\frac{1-\xi}\xi>1$ for $\xi\in (\frac 13,\frac 12)$. Therefore, the fixed point $v_1^*$ is still unstable, but it does not oscillate since $(\frac{1-\xi}\xi)^2 >0$. Actually we can find two fixed points for \refe{dm2squared_pm}:
\bi
\item If $v_1^*$ is decreased; i.e., if $\tilde v_1(t)<0$, then from \refe{dm2squared_pm} we can see that $v_3(t+T)$ increases and $v_1(t+2T)$ decreases. However, since $v_3(t+T)$ cannot be greater than 1, \refe{dm2squared_pm} has a new fixed point: $v_3^*=1$ and $v_1^*=2-\frac{1-\xi} \xi$, which is between 0 and 1.
\item Similarly, if $v_1^*$ is increased; i.e., if $\tilde v_1(t)>0$, \refe{dm2squared_pm} has another fixed point: $v_1^*=1$ and $v_3^*=2-\frac{1-\xi} \xi$.
\ei
These two fixed points are stable. These results can be easily confirmed by CTM simulations.

\begin{figure} \bc
$\ba{c@{\hspace{0.3in}}c}
\includegraphics[height=2.5in]{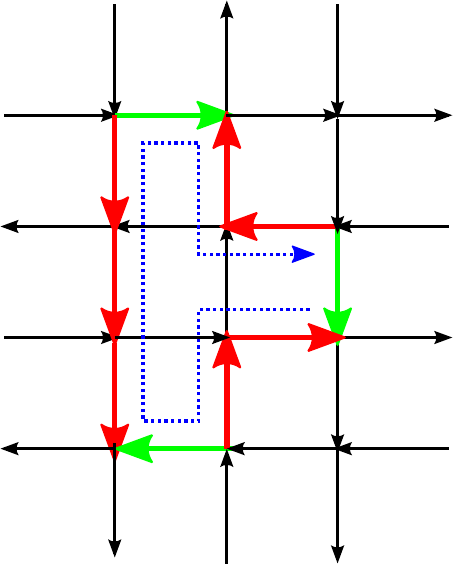} &
\includegraphics[height=2.5in]{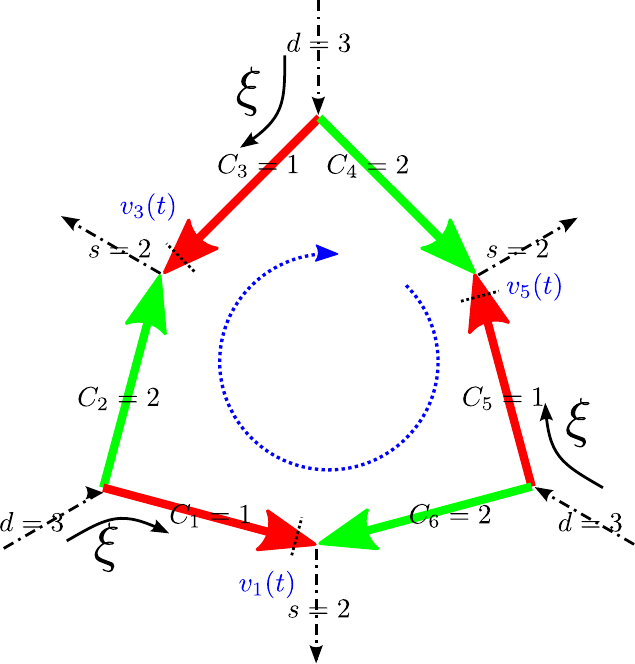} \\
\multicolumn{1}{c}{\mbox{\bf (a)}} &
    \multicolumn{1}{c}{\mbox{\bf (b)}}
\ea$
\caption{A $4\times 3$ grid network (a) and the $(DM)^3$ network (b)} \label{4by3dm2cubed} \ec 
\end{figure}

Third, we consider a $4\times 3$ grid network in \reff{4by3dm2cubed}(a), which embeds a $(DM)^3$ network shown in \reff{4by3dm2cubed}(b). Under certain conditions, it is possible that there exists a circular information propagation path formed by congested (red) and uncongested (green) links. For the $(DM)^3$ network, we have the same set-up as in the $(DM)^2$ network in \reff{2by2dm2squared}(b). For $\xi\in(\frac 13,\frac 12)$, we can then obtain the following \Po map:
\bsq \label{dm2cubed_pm}
\bqn
v_1(t+T)&=&\min\{1,2-\frac{1-\xi}\xi v_5(t)\},\\
v_3(t+T)&=&\min\{1,2-\frac{1-\xi}\xi v_1(t)\},\\
v_5(t+T)&=&\min\{1,2-\frac{1-\xi}\xi v_3(t)\},
\eqn
\esq
which has a fixed point $v_1^*=v_3^*=v_5^*=2\xi$. For a small perturbation $\tilde v_1(t)$ around the fixed point, the corresponding map becomes
\bqs
\tilde v_1(t+3T)&=&-(\frac{1-\xi}\xi)^3 \tilde v_1(t),
\eqs
which is unstable, since $\frac{1-\xi}\xi>1$ for $\xi\in (\frac 13,\frac 12)$, and oscillates, since $-(\frac{1-\xi}\xi)^3<0$. Therefore PPO traffic patterns can develop in this network.

From the analyses above, we can make the following conclusions regarding a general  $(DM)^n$ ($n\geq 1$) network with a symmetric set-up as in \reff{2by2dm2squared}(b) and \reff{4by3dm2cubed}(b) and $\xi\in(\frac 13,\frac12)$: (i) The network can be stationary at $v_1^*=2\xi$, which is unstable, and the perturbation $\tilde v_1(t)$ has the following dynamics:
\bqn
\tilde v_1(t+nT)&=&-(\frac{1-\xi}\xi)^n \tilde v_1(t).
\eqn
(ii) When $n$ is odd, asymptotic PPO traffic patterns can appear under any non-stationary initial conditions. (iii) When $n$ is even, asymptotic PPO traffic patterns can appear under symmetric non-stationary initial conditions, but two more stable stationary states can appear under asymmetric non-stationary initial conditions. (iv) These networks can be embedded in a large grid network. (v) All these observations can be confirmed with CTM simulations.

\subsection{Beltway networks}
\begin{figure} \bc
$\ba{c@{\hspace{0.3in}}c}
\includegraphics[height=2.5in]{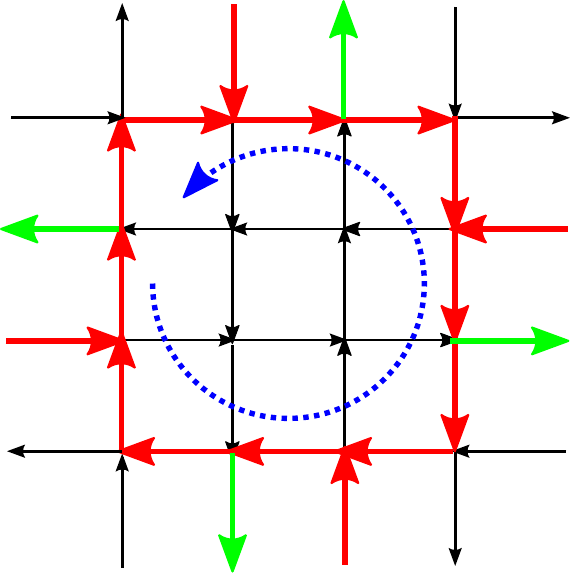} &
\includegraphics[height=2.5in]{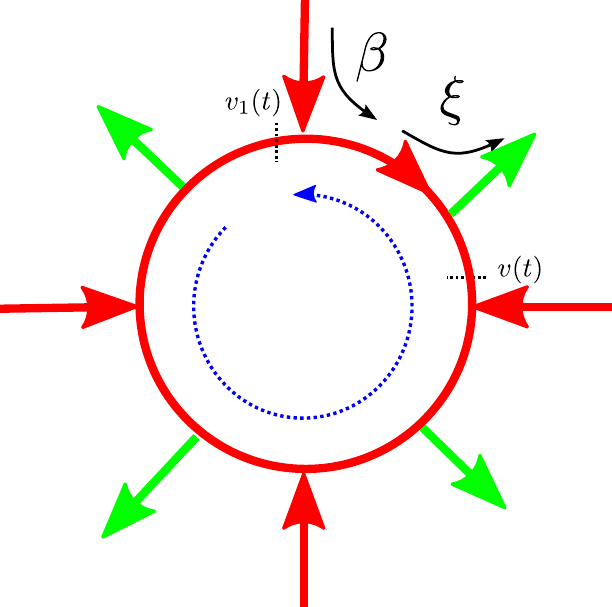} \\
\multicolumn{1}{c}{\mbox{\bf (a)}} &
    \multicolumn{1}{c}{\mbox{\bf (b)}}
\ea$
\caption{A $4\times 4$ grid network (a) and a beltway network (b)} \label{4by4beltway} \ec 
\end{figure}

In the $(DM)^n$ networks, circular information propagation paths consist of both congested and uncongested links. In this subsection, we consider another type of circular information propagation paths consisting of all congested links on a ring road. A congested ring road with on- and off-ramps is called a beltway, whose traffic dynamics were first studied in \citep{daganzo1996gridlock}. It was revealed that the beltway network can become totally gridlocked. In \citep{daganzo2007gridlock}, the gridlock mechanism was also discussed for urban networks. In this subsection we apply the \Po map approach to analyze traffic dynamics in a beltway network, shown in \reff{4by4beltway}(b), which can be embedded in a grid network. Here we consider a symmetric beltway network with $n$ pairs of alternate on- and off-ramps, in which the ring road and on-ramps are congested, and the off-ramps are not. In addition, we apply the priority-based merge model in \refe{mergemodel}, where the on-ramp's merging ratio is $\beta$, and the FIFO diverge model in \refe{divergemodel}, where the turning proportion to the off-ramp is $\xi$.

As shown in \reff{4by4beltway}, we define two \Po sections at the upstream points of two consecutive merges and denote the two out-fluxes of the mainline road by $v(t)$ and $v_1(t)$ respectively. Then after a pair of off- and on-ramps, the out-flux becomes $v_1(t+T)=\frac{1-\beta}{1-\xi}v(t)$, where $T$ is the time for the shock wave to travel . After $n$ pairs of off- and on-ramps, we obtain the following \Po map:
\bqn
v(t+nT)&=&(\frac{1-\beta}{1-\xi})^n v(t). \label{beltway_pm}
\eqn
In \citep{daganzo1996gridlock}, two different parameters are used: $\alpha=\frac\beta{1-\beta}$, and $\mu=\frac{\xi}{1-\xi}$. Then the coefficient $\frac{1-\beta}{1-\xi}=\frac{1+\mu}{1+\alpha}$, which is given in equation (2) of \citep{daganzo1996gridlock}. Note that the mapping in \citep{daganzo1996gridlock} was implicitly derived ``once the observer has traveled around the loop once''. This is different from our approach based on circular information propagation.

From the \Po map, \refe{beltway_pm}, we have the following observations: (i) When $\frac{1-\beta}{1-\xi}<1$, the network converges to the gridlock state ($v(t)=0$), and the gridlock state is stable. (ii) When $\frac{1-\beta}{1-\xi}>1$, the gridlock state is still a stationary state, but it is unstable. (iii) When $\frac{1-\beta}{1-\xi}=1$, there can exist multiple stationary states. Furthermore, other results, including the flow half-life and flow recovery, in \citep{daganzo1996gridlock} can also be obtained (but omitted here) from the \Po map, \refe{beltway_pm}.

\begin{figure} \bc
$\ba{c@{\hspace{0.3in}}c}
\includegraphics[height=2.5in]{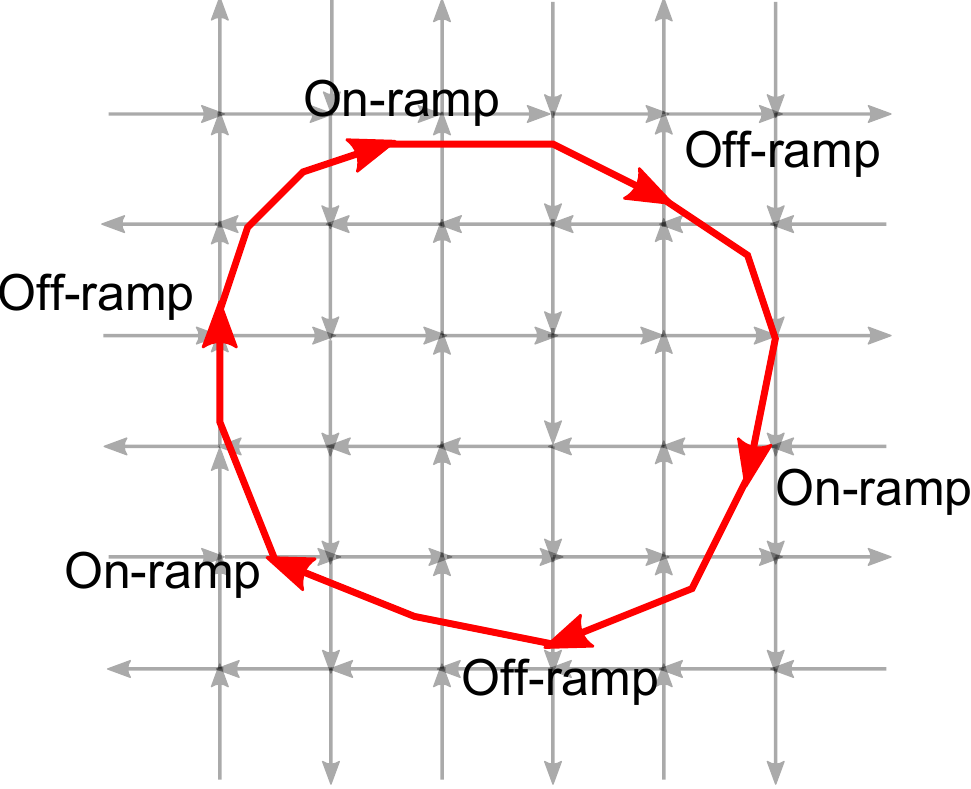} &
\includegraphics[height=2.5in]{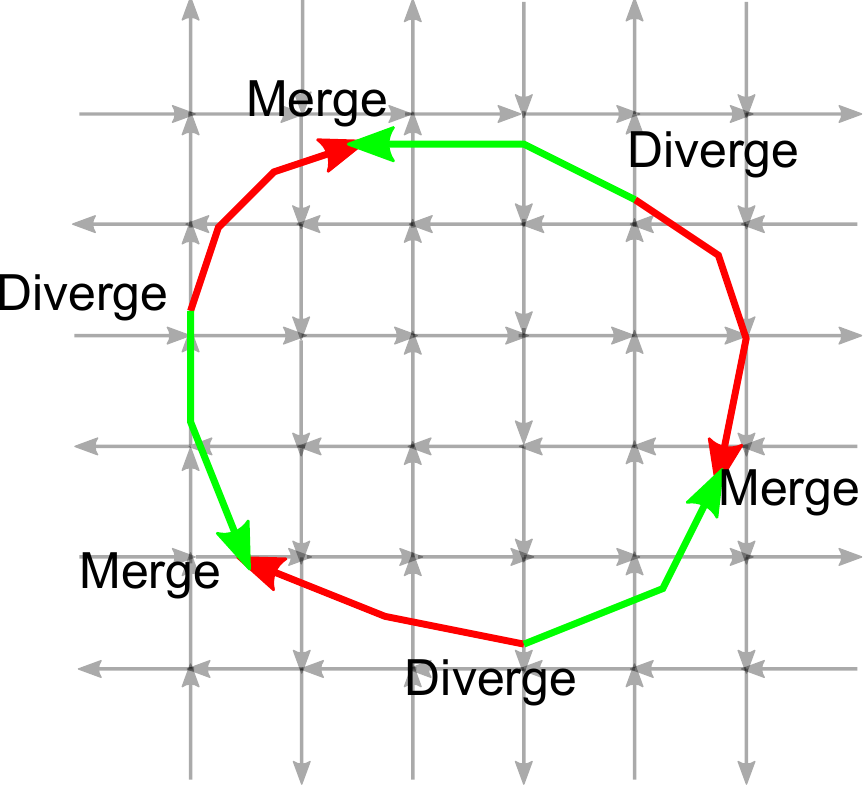} \\
\multicolumn{1}{c}{\mbox{\bf (a)}} &
    \multicolumn{1}{c}{\mbox{\bf (b)}}
\ea$
\caption{Two types of networks with circular information propagation: beltway network (a) and $(DM)^n$ network (b)} \label{circularprop} \ec 
\end{figure}

Moreover, we can see that both network-induced unstable traffic and gridlock have to be associated with circular information propagation. As shown in \reff{circularprop}, there can only be two types of (nontrivial) circular information propagation in a road network: all links are congested and form a beltway network, or links are alternatively congested and uncongested and form a $(DM)^n$ network. \footnote{Of course when all links on a ring are uncongested, there can also exists a circular path for information propagation, but this case is trivial and not interesting.} Therefore, the $(DM)^n$ network structure is sufficient and necessary conditions for the existence of network-induced unstable traffic patterns, and the beltway network structure is sufficient and necessary conditions for the existence of network-induced gridlock states.

\section{Conclusion}

In this study, we first studied the global traffic dynamics arising in the kinematic wave model for a diverge-merge (DM) network with two intermediate links when the downstream link constitutes a bottleneck. In stationary states, the two intermediate links admit respectively under-critical and over-critical stationary states, in which the kinematic waves propagate in a circular fashion. Then we derived a \Po map for flow-rates at two \Po sections. For the \Po map, which is a one-dimensional discrete dynamical system, we studied the fixed points and their stability properties. We showed that the fixed points can be unstable under the conditions when one intermediate link is strictly under-critical (SUC) and the other is strictly over-critical (SOC). Furthermore, we found that the \Po map can have  periodical points of period 2 when it is unstable, but it has no other periodical points or chaos. We also demonstrated that the route choice proportion can cause a Hopf bifurcation in the stability of fixed points. By comparing the results for the \Po map with those for the kinematic wave model in previous studies, we found that the \Po map indeed describes the global traffic dynamics in the DM network: the fixed points of the \Po map correspond to stationary link flow-rates, their stability is related to the stability of stationary states, and the periodical points of period 2 yield the oscillation magnitudes in the flow-rates for persistent periodic oscillatory traffic patterns.
Thus together with \citep{jin2012statics}, this study provides us a more complete picture regarding traffic dynamics arising in  the kinematic wave model of traffic dynamics in the DM network: stationary states always exist under constant boundary conditions, but a stationary state can be finite-time stable, asymptotically stable, or unstable; when it is unstable, the oscillation magnitudes are still bounded, and chaos does not arise in the DM network.
We further applied the \Po map approach to analyzing traffic patterns in more general $(DM)^n$ and beltway networks, which can be embedded in grid networks, since a \Po map can still be derived even when a circular path involves many links and junctions. 

Through this study, we obtain the following insights into network traffic dynamics as well as the \Po map approach: (i) The $(DM)^n$ and beltway networks are the only networks that admit nontrivial circular information propagation; (ii) The Poincare map approach can be applied on circular information propagation paths with any number of junctions and links; (iii) The $(DM)^n$ networks are sufficient and necessary structures for network-induced unstable traffic patterns; and (iv) The beltway networks are sufficient and necessary structures for network-induced gridlock states. 
Some of the insights on network traffic flow have been obtained in the literature with the kinematic wave model in \citep{daganzo1996gridlock,jin2009network}, but the \Po map approach is much simpler and more powerful, as it can be applied for very general networks and lead to much broader and deeper insights on network traffic flow.

In the future, we will be interested in collecting empirical evidences for the existence of such network-induced instability. 
From the analysis we can see that the stability property of traffic dynamics in a road network is related to the network structure, link capacities, route choice proportions, and merging ratios. Therefore, incidents, route guidance strategies, and traffic signals can all cause unstable traffic patterns. In the future, we will be interested in various controlling strategies to stabilize traffic in a road network.  
Another research topic is to examine the stability of traffic dynamics in a road network subject to both route choice behaviors and interactions among network bottlenecks, since unstable traffic dynamics could have impacts on users' route choice behaviors and user equilibrium in the network.

\section*{Acknowledgments}
We would like to thank two anonymous reviewers for their constructive comments. The views and results are the author's alone.

\section*{Appendix A. Some basic properties of min and max operators}
The following are some basic properties of min and max operators. Here $a$, $b$, $c$, and $\la$ are all real numbers.
\bqs
\min\{a,b\}&=&-\max\{-a,-b\}\\
\max\{a,b\}&=&-\min\{-a,-b\}\\
\min\{a,\max\{b,c\}\}&=&\max\{\min\{a,b\},\min\{a,c\}\}\\
\max\{a,\min\{b,c\}\}&=&\min\{\max\{a,b\},\max\{a,c\}\}\\
\min\{a,\min\{b,c\}\}&=&\min\{a,b,c\}\\
\max\{a,\max\{b,c\}\}&=&\max\{a,b,c\}\\
a+\min\{b,c\}&=&\min\{a+b,a+c\}\\
a+\max\{b,c\}&=&\max\{a+b,a+c\}\\
\la \min\{a,b\}&=&\min\{\la a, \la b\}, \la>0\\
\la \max\{a,b\}&=&\max\{\la a, \la b\}, \la>0
\eqs

\section*{Appendix B. Proof of Theorem \ref{thm:fts}}
{\em Proof}. We first prove the results for $\xi\in \Xi_1$.
\ben
 \item When $C_3\leq C_0$, $C_3<C_1+C_2$, and $\frac {C_1}{C_3}\leq \xi\leq 1$. For any initial $v\leq C_1 \leq \xi C_3$, we have $C_3-\frac{1-\xi}{\xi} v \geq \xi C_3 \geq C_1$. Thus from \refe{pmap} we have $F v=C_1$, which is the fixed point. Thus the \Po map converges to the fixed point in one time step.

\item When $C_3\leq C_0$, $C_3<C_1+C_2$, $\xi\in(1-\frac{C_2}{C_3}, \frac{C_1}{C_3})$, and $\xi=\beta$. Since $C_3-C_2<\xi C_3$ and $C_3-(1-\xi) C_0\leq \xi C_3$, we have
\bqs
F v&=&\min\{C_1,\max\{\xi C_3, C_3-\frac{1-\xi}{\xi} v\}\}.
\eqs
If initially $v\geq \xi C_3$, then $C_3-\frac{1-\xi}{\xi} v\leq \xi C_3$, and  $Fv=\xi C_3 =v^* <C_1$. If initially $v< \xi C_3$, then $C_3-\frac{1-\xi}{\xi} v> \xi C_3$, $F v=\min\{C_1, C_3-\frac{1-\xi}{\xi} v\}>\xi C_3$, and $F^2 v=\xi C_3=v^*$. Thus the \Po map converges to the fixed point after at most two time-steps. 

\item When $C_3=C_0<C_1+C_2$, $\xi\in(1-\frac{C_2}{C_3}, \frac{C_1}{C_3})$, and $\xi>\beta$. Since $C_3-C_2<\xi C_3$ and $\beta C_3<\xi C_3$, we have
\bqs
F v&=&\min\{C_1,\max\{\xi C_3, C_3-\frac{1-\xi}{\xi} v\}\}.
\eqs
If initially  $v\geq \xi C_3$, then $C_3-\frac{1-\xi}{\xi} v\leq \xi C_3$, and  $Fv=\xi C_3 =v^*<C_1$. If initially $v< \xi C_3$, then $C_3-\frac{1-\xi}{\xi} v> \xi C_3$, $F v=\min\{C_1, C_3-\frac{1-\xi}{\xi} v\}>\xi C_3$, and $F^2 v=\xi C_3=v^*$. Thus the \Po map converges to the fixed point after at most two time-steps.  
\een
Therefore, the fixed points of \refe{pmap} are finite-time stable for $\xi\in \Xi_1$. Similarly we can prove the results for $\xi\in\Xi_2$. \eop

\section*{Appendix C. Proofs of Lemmas in Section 4}
In the following proofs, we denote $\la=\frac{1-\xi}{\xi}$.

\bi
\item {\em Proof of Lemma \ref{lemma:p1}}. 

Since $\xi>\frac 12$, $\la<1$. 

Assume that \refe{eqn:period2} has a solution $v<\xi C_3$, then $C_3-\la v>\xi C_3>A_1$. From \refe{eqn:period2} we have
\bqs
v&=&\min\{C_1, \max\{A_1, C_3-\la\min\{C_1, C_3-\la v\}\}\}\\&=&\min\{C_1,\max\{A_1, C_3-\la C_1, C_3-\la(C_3-\la v)\}\}.
\eqs 
Since $A_1<C_1$, $C_3-\la C_1<C_1$, and $C_3-\la (C_3-\la v)<\xi C_3<C_1$, we have
\bqs
v&=&\max\{A_1, C_3-\la C_1, C_3-\la(C_3-\la v)\}.
\eqs
If $v=C_3-\la(C_3-\la v)$, then $v=\xi C_3$. Since $v<\xi C_3$, the only possible solution is $v=\max\{A_1,C_3-\la C_1\}<\xi C_3$. But for this solution to exist, we require that $C_3-\la(C_3-\la v)\leq \max\{A_1,C_3-\la C_1\}=v$,
which leads to $(1-\la) (C_3-v(1+\la))=(1-\la)(C_3-\frac {v} {\xi})\leq 0$.
Since $v<\xi C_3$, we require that $1-\la\leq 0$, which is equivalent to that $\la\geq 1$ or $\xi\leq \frac 12$. This contradicts that $\xi>\frac 12$.

Similarly, we can show that  \refe{eqn:period2} has no solution of $v>\xi C_3$.
\eop

\item {\em Proof of Lemma \ref{lemma:p2}}. 

Since $\xi\leq \frac 12$, $\la\geq 1$. 

When $C_3<\min\{C_0,C_1+C_2\}$, $\xi\in\tilde\Xi_1$, and $\xi< \frac 12$, it is straightforward to verify that $v^-<v^*<v^+$. In addition, 
\bqs
F v^-&=&\min\{C_1, \max\{A_1,C_3-\la v^-\}\}=\min\{C_1, C_3-\la v^-\}\\&=&\min\{C_1, C_3-\la A_1, (1-\la)C_3+\la^2 C_1\}=v^+,
\eqs
since $C_1\leq (1-\la)C_3+\la^2 C_1$ for $\la\geq 1$.
Similarly we can show that $F v^+=v^-$. Thus $v^-$ and $v^+$ are two periodical points of period 2. 

Next we show that they are the only periodical points of period 2 when $\xi<\frac 12$.
First, \refe{eqn:period2} can be simplified as follows:
\bqs
v&=&\min\{C_1,\max\{A_1,C_3-\la C_1, C_3-\la \max\{A_1,C_3-\la v\}\}\}\\
&=&\max\{\min\{C_1,v^-\},\min\{C_1,C_3-\la A_1, C_3-\la(C_3-\la v)\}\}\\
&=&\max\{\min\{C_1,v^-\},\min\{v^+,(1-\la)C_3+\la^2 v\}\}.
\eqs
Since $A_1<C_1$ and $C_3-\la C_1<C_1$, \refe{eqn:period2} is equivalent to
\bqn
v&=&\max\{v^-, \min\{v^+, (1-\la)C_3+\la^2 v\}\}. \label{eqn:period2s}
\eqn
When $\xi<\frac 12$; i.e., when $\la>1$, we have that $v>,=,<(1-\la)C_3+\la^2 v$ if and only if $v<,=,>v^*$, respectively. If $v<v^*$, then $(1-\la)C_3+\la^2 v<v<v^*<v^+$, and \refe{eqn:period2s} leads to $v=\max\{v^-, (1-\la)C_3+\la^2 v\}=v^-$; If $v>v^*$, then 
$(1-\la)C_3+\la^2 v>v>v^*>v^-$, and \refe{eqn:period2s} leads to $v=\min\{v^+, (1-\la)C_3+\la^2 v\}=v^+$. Thus $v^-$ and $v^+$ are the only two periodical points of period 2 when $\xi<\frac 12$.

When $\xi=\frac 12$; i.e., when $\la=1$, \refe{eqn:period2s} is equivalent to $v=\max\{v^-,\min\{v^+,v\}\}$. In this case, $v^-=\max\{A_1,C_3- C_1\}$, and $v^+=\min\{C_1, C_3-A_1\}$. Thus, any $v\in[v^-,v^*)\cup (v^*,v^+]$ is a periodical point of period 2.
\eop

\item {\em Proof of Lemma \ref{lemma:p3}}. 

Since $\xi\leq \frac 12$, $\la\geq 1$. 

Assume that there exists a $v<\xi C_3$ such that $F^4v =v$. Then we have 
\bqs
F^2v&=&\max\{A_1, C_3-\la Fv\}=\max\{A_1, C_3-\la C_1, C_3-\la(C_3-\la v)\},\\
F^4v&=&\max\{A_1, C_3-\la C_1, C_3-\la(C_3-\la F^2 v)\}\\
&=&\max\{A_1, C_3-\la C_1, C_3(1-\la)+\la^2 A_1, C_3(1-\la)+\la^2 (C_3-\la C_1),\\
&&C_3(1-\la)+\la^2(C_3-\la(C_3-\la v))\}.
\eqs
Since $F^4 v=v$,  $v=C_3(1-\la)+\la^2(C_3-\la(C_3-\la v))$, which leads to $v=\xi C_3$, or $v=\max\{A_1, C_3-\la C_1, C_3(1-\la)+\la^2 A_1, C_3(1-\la)+\la^2 (C_3-\la C_1)\}$. 
Since $v<\xi C_3$, the only possible solution is 
\bqs
v&=&\max\{A_1, C_3-\la C_1, C_3(1-\la)+\la^2 A_1, C_3(1-\la)+\la^2 (C_3-\la C_1)\}.
\eqs
Since $\la \geq 1$, $A_1<\xi C_3$, and $\xi C_3<C_1$, we have $A_1\geq C_3(1-\la)+\la^2 A_1$ and $C_3-\la C_1 \geq C_3(1-\la)+\la^2 (C_3-\la C_1)$. Therefore we have
\bqs
v&=&\max\{A_1,C_3-\la C_1\}=v^-,
\eqs
which is a periodical point of period 2.

Similarly, if $v>\xi C_3$ such that $F^4 v=v$, then $v=v^+$, which is also a periodical point of period 2.
In both cases, the solutions of $F^4 v=v$ are also those of $F^2v=v$. Thus, there are no periodical points of period 4 for the \Po  map, \refe{pmap}. 
 \eop

\ei


\bibliographystyle{elsarticle-harv}

\begin{thebibliography}{33}
\expandafter\ifx\csname natexlab\endcsname\relax\def\natexlab#1{#1}\fi
\expandafter\ifx\csname url\endcsname\relax
  \def\url#1{\texttt{#1}}\fi
\expandafter\ifx\csname urlprefix\endcsname\relax\def\urlprefix{URL }\fi

\bibitem[{Ahn and Cassidy(2007)}]{ahn2007freeway}
Ahn, S., Cassidy, M., 2007. Freeway traffic oscillations and vehicle
  lane-change maneuvers. In: Proceedings of the 17th International Symposium on
  Traffic and Transportation Theory. pp. 691--710.

\bibitem[{Ahn et~al.(2010)Ahn, Laval, and Cassidy}]{ahn2010merging}
Ahn, S., Laval, J., Cassidy, M.~J., 2010. Merging and diverging effects on
  freeway traffic oscillations: Theory and observation. Transportation Research
  Record: Journal of the Transportation Research Board 2188, 1--8.

\bibitem[{Barth and Boriboonsomsin(2008)}]{barth2008real}
Barth, M., Boriboonsomsin, K., 2008. {Real-world carbon dioxide impacts of
  traffic congestion}. Transportation Research Record: Journal of the
  Transportation Research Board 2058, 163--171.

\bibitem[{Bhat and Bernstein(2000)}]{bhat2000finite}
Bhat, S., Bernstein, D., 2000. {Finite-time stability of continuous autonomous
  systems}. SIAM Journal on Control and Optimization 38~(3), 751--766.

\bibitem[{Bressan(1996)}]{bressan1996semigroup}
Bressan, A., 1996. {The semigroup approach to systems of conservation laws}.
  Matematica Contemporanea 10, 21--74.

\bibitem[{Burns and Hasselblatt(2011)}]{burns2011sharkovsky}
Burns, K., Hasselblatt, B., 2011. {The Sharkovsky theorem: A natural direct
  proof}. American Mathematical Monthly 118~(3), 229--244.

\bibitem[{Daganzo(2007)}]{daganzo2007gridlock}
Daganzo, C., 2007. {Urban gridlock: Macroscopic modeling and mitigation
  approaches}. Transportation Research Part B 41~(1), 49--62.

\bibitem[{Daganzo(1995)}]{daganzo1995ctm}
Daganzo, C.~F., 1995. {The cell transmission model {II}: Network traffic}.
  Transportation Research Part B 29~(2), 79--93.

\bibitem[{Daganzo(1996)}]{daganzo1996gridlock}
Daganzo, C.~F., 1996. The nature of freeway gridlock and how to prevent it.
  Proceedings of the 13th International Symposium on Transportation and Traffic
  Theory, 629--646.

\bibitem[{Daganzo(1999)}]{daganzo1999remarks}
Daganzo, C.~F., 1999. {Remarks on traffic flow modeling and its applications}.
  In: Brilon, W., Huber, F., Schreckenberg, M., Wallentowitz, H. (Eds.),
  Proceedings of Traffic and Mobility: Simulation, Economics and Environment.
  Springer Verlag, pp. 105--115.

\bibitem[{Daganzo et~al.(1999)Daganzo, Cassidy, and Bertini}]{daganzo1999phase}
Daganzo, C.~F., Cassidy, M.~J., Bertini, R.~L., 1999. Possible explanations of
  phase transitions in highway traffic. Transportation Research A 33~(5),
  365--379.

\bibitem[{Daganzo et~al.(2011)Daganzo, Gayah, and
  Gonzales}]{daganzo2011bifurcations}
Daganzo, C.~F., Gayah, V.~V., Gonzales, E.~J., 2011. Macroscopic relations of
  urban traffic variables: Bifurcations, multivaluedness and instability.
  Transportation Research Part B 45~(1), 278--288.

\bibitem[{Engquist and Osher(1980)}]{engquist1980difference}
Engquist, B., Osher, S., 1980. {One-sided difference schemes and transonic
  flow}. Proceedings of the National Academy of Sciences 77~(6), 3071--3074.

\bibitem[{Galor(2004)}]{galor2004introduction}
Galor, O., 2004. {Introduction to Stability Analysis of Discrete Dynamical
  Systems}. Macroeconomics.

\bibitem[{Greenshields(1935)}]{greenshields1935capacity}
Greenshields, B.~D., 1935. A study in highway capacity. Highway Research Board
  Proceedings 14, 448--477.

\bibitem[{Haimo(1986)}]{haimo1986finite}
Haimo, V., 1986. {Finite time controllers}. SIAM Journal on Control and
  Optimization 24, 760.

\bibitem[{Herman et~al.(1959)Herman, Montroll, Potts, and
  Rothery}]{herman1959traffic}
Herman, R., Montroll, E., Potts, R., Rothery, R., 1959. {Traffic dynamics:
  analysis of stability in car following}. Operations research 7~(1), 86--106.

\bibitem[{Holmgren(1996)}]{holmgren1996first}
Holmgren, R., 1996. {A first course in discrete dynamical systems}. Springer
  Verlag.

\bibitem[{Jin(2003)}]{jin2003dissertation}
Jin, W.-L., 2003. Kinematic wave models of network vehicular traffic. Ph.D.
  thesis, University of California, Davis.
\newline\urlprefix\url{http://arxiv.org/abs/math.DS/0309060}

\bibitem[{Jin(2009)}]{jin2009network}
Jin, W.-L., 2009. {Asymptotic traffic dynamics arising in diverge-merge
  networks with two intermediate links}. Transportation Research Part B 43~(5),
  575--595.

\bibitem[{Jin(2012{\natexlab{a}})}]{jin2012network}
Jin, W.-L., 2012{\natexlab{a}}. A kinematic wave theory of multi-commodity
  network traffic flow. Transportation Research Part B 46~(8), 1000--1022.

\bibitem[{Jin(2012{\natexlab{b}})}]{jin2012statics}
Jin, W.-L., 2012{\natexlab{b}}. The traffic statics problem in a road network.
  Transportation Research Part B 46~(10), 1360–1373.

\bibitem[{Jin and Zhang(2005)}]{jin2005paramics}
Jin, W.-L., Zhang, Y., 2005. Paramics simulation of periodic oscillations
  caused by network geometry. Transportation Research Record: Journal of the
  Transportation Research Board 1934, 188--196.

\bibitem[{Kerner and Konh\"auser(1994)}]{kerner1994cluster}
Kerner, B.~S., Konh\"auser, P., 1994. Structure and parameters of clusters in
  traffic flow. Physical Review E 50~(1), 54--83.

\bibitem[{LaSalle(1976)}]{lasalle1976stability}
LaSalle, J., 1976. {The stability of dynamical systems}. Society for Industrial
  Mathematics.

\bibitem[{Lebacque(1996)}]{lebacque1996godunov}
Lebacque, J.~P., 1996. {The Godunov scheme and what it means for first order
  traffic flow models}. Proceedings of the 13th International Symposium on
  Transportation and Traffic Theory, 647--678.

\bibitem[{Li and Yorke(1975)}]{li1975period}
Li, T., Yorke, J., 1975. {Period three implies chaos}. American mathematical
  monthly, 985--992.

\bibitem[{Li et~al.(2010)Li, Peng, and Ouyang}]{li2010measurement}
Li, X., Peng, F., Ouyang, Y., 2010. {Measurement and estimation of traffic
  oscillation properties}. Transportation Research Part B 44~(1), 1--14.

\bibitem[{Mauch and Cassidy(2002)}]{mauch2002freeway}
Mauch, M., Cassidy, M.~J., 2002. Freeway traffic oscillations: observations and
  predictions. In: Proceedings of the 15th International Symposium on
  Transportation and Traffic Theory.

\bibitem[{Oh et~al.(2001)Oh, Oh, Ritchie, and Chang}]{oh2001real}
Oh, C., Oh, J.~S., Ritchie, S.~G., Chang, M., 2001. Real-time estimation of
  freeway accident likelihood. In: Proceedings of the 80th Annual Meeting of
  the Transportation Research Board.

\bibitem[{Slotine and Li(1991)}]{slotine1991applied}
Slotine, J., Li, W., 1991. {Applied nonlinear control}. Vol. 461. Prentice hall
  Englewood Cliffs, NJ.

\bibitem[{Sugiyama et~al.(2008)Sugiyama, Fukui, Kikuchi, Hasebe, Nakayama,
  Nishinari, Tadaki, and Yukawa}]{sugiyama2008traffic}
Sugiyama, Y., Fukui, M., Kikuchi, M., Hasebe, K., Nakayama, A., Nishinari, K.,
  Tadaki, S., Yukawa, S., 2008. {Traffic jams without
  bottlenecks�experimental evidence for the physical mechanism of the
  formation of a jam}. New Journal of Physics 10, 033001.

\bibitem[{Wiggins and Heck(2003)}]{wiggins2003chaos}
Wiggins, S., Heck, A., 2003. {Introduction to Applied Nonlinear Dynamical
  Systems and Chaos}. Springer.

\end{thebibliography}

\end {document}